%%% -*- plain-TeX -*-
%%% V. Ancona,  M. Maggesi
%%% On the quantum cohomology of Fano bundles over projective spaces

%% @texfile{
%%   author = "Karl Berry",
%%   version = "2.8-beta",
%%   date = "Sun May  7 07:32:53 EDT 2000",
%%   filename = "eplain.tex",
%%   email = "karl@cs.umb.edu",
%%   address = "135 Center Hill Rd. // Plymouth, MA 02360"
%%   checksum = "2369    5896   80103",
%%   codetable = "ISO/ASCII",
%%   supported = "yes",
%%   docstring = "This file defines macros that extend and expand on
%%                plain TeX. eplain.tex is xeplain.tex and the other
%%                source files with comments stripped; see the original
%%                files for author credits, etc.  And please base diffs
%%                or other contributions on xeplain.tex, not the
%%                stripped-down eplain.tex.",
%% }
\ifx\eplain\undefined
  \let\next\relax
\else
  \expandafter\let\expandafter\next\csname endinput\endcsname
\fi
\next
\def\makeactive#1{\catcode`#1 = \active \ignorespaces}%
\chardef\letter = 11
\chardef\other = 12
\edef\leftdisplays{\the\catcode`@}%
\catcode`@ = \letter
\let\@eplainoldatcode = \leftdisplays
\toksdef\toks@ii = 2
\def\uncatcodespecials{%
   \def\do##1{\catcode`##1 = \other}%
   \dospecials
}%
{%
   \makeactive\^^M %
   \long\gdef\letreturn#1{\let^^M = #1}%
}%
\let\@eattoken = \relax  % Define this, so \eattoken can be used in \edef.
\def\eattoken{\let\@eattoken = }%
\def\gobble#1{}%
\def\gobbletwo#1#2{}%
\def\gobblethree#1#2#3{}%
\def\@emptymarkA{\@emptymarkB}
\def\ifempty#1{\@@ifempty #1\@emptymarkA\@emptymarkB}%
\def\@@ifempty#1#2\@emptymarkB{\ifx #1\@emptymarkA}%
\def\@gobblemeaning#1:->{}%
\def\sanitize{\expandafter\@gobblemeaning\meaning}%
\def\ifundefined#1{\expandafter\ifx\csname#1\endcsname\relax}%
\def\csn#1{\csname#1\endcsname}%
\def\ece#1#2{\expandafter#1\csname#2\endcsname}%
\def\expandonce{\expandafter\noexpand}%
\let\@plainwlog = \wlog
\let\wlog = \gobble
\newlinechar = `^^J
\def\loggingall{\tracingcommands\tw@\tracingstats\tw@
   \tracingpages\@ne\tracingoutput\@ne\tracinglostchars\@ne
   \tracingmacros\tw@\tracingparagraphs\@ne\tracingrestores\@ne
   \showboxbreadth\maxdimen\showboxdepth\maxdimen
}%
\def\tracingoff{\tracingonline\z@\tracingcommands\z@\tracingstats\z@
  \tracingpages\z@\tracingoutput\z@\tracinglostchars\z@
  \tracingmacros\z@\tracingparagraphs\z@\tracingrestores\z@
  \showboxbreadth5 \showboxdepth3
}%
\begingroup
  \catcode`\{ = 12 \catcode`\} = 12
  \catcode`\[ = 1 \catcode`\] = 2
  \gdef\lbracechar[{]%
  \gdef\rbracechar[}]%
  \catcode`\% = \other
  \gdef\percentchar[%]\endgroup
\def^^L{\par}%
\def\vpenalty{\ifhmode\par\fi \penalty}%
\def\hpenalty{\ifvmode\leavevmode\fi \penalty}%
\def\iterate{%
  \let\loop@next\relax
  \body
  \let\loop@next\iterate
  \fi
  \loop@next
}%
\def\edefappend#1#2{%
  \toks@ = \expandafter{#1}%
  \edef#1{\the\toks@ #2}%
}%
\long\def\hookprepend{\@hookassign{\the\toks@ii \the\toks@}}%
\long\def\hookappend{\@hookassign{\the\toks@ \the\toks@ii}}%
\let\hookaction = \hookappend % either one should be ok
\long\def\@hookassign#1#2#3{%
  \expandafter\ifx\csname @#2hook\endcsname \relax
    \toks@ = {}%
  \else
    \expandafter\let\expandafter\temp \csname @#2hook\endcsname
    \toks@ = \expandafter{\temp}%
  \fi
  \toks2 = {#3}% Don't expand the argument all the way.
  \ece\edef{@#2hook}{#1}%
}%
\long\def\hookactiononce#1#2{%
  \edefappend#2{\global\let\noexpand#2\relax}
  \hookaction{#1}#2%
}%
\def\hookrun#1{%
  \expandafter\ifx\csname @#1hook\endcsname \relax \else
    \def\temp{\csname @#1hook\endcsname}%
    \expandafter\temp
  \fi
}%
\def\setproperty#1#2#3{\ece\edef{#1@p#2}{#3}}%
\def\getproperty#1#2{%
  \expandafter\ifx\csname#1@p#2\endcsname\relax
  \else \csname#1@p#2\endcsname
  \fi
}%
\ifx\@undefinedmessage\@undefined
  \def\@undefinedmessage
    {No .aux file; I won't warn you about undefined labels.}%
\fi
%% @texfile{
%%   author = "Karl Berry and Oren Patashnik",
%%   version = "0.99k",
%%   date = "13 Nov 1995",
%%   filename = "btxmac.tex",
%%   address = "Please use electronic mail",
%%   checksum = "842    4571   33524",
%%   email = "opbibtex@cs.stanford.edu",
%%   codetable = "ISO/ASCII",
%%   supported = "yes",
%%   docstring = "Defines macros that make BibTeX work with plain TeX",
%% }
\edef\cite{\the\catcode`@}%
\catcode`@ = 11
\let\@oldatcatcode = \cite
\chardef\@letter = 11
\chardef\@other = 12
\def\@innerdef#1#2{\edef#1{\expandafter\noexpand\csname #2\endcsname}}%
\@innerdef\@innernewcount{newcount}%
\@innerdef\@innernewdimen{newdimen}%
\@innerdef\@innernewif{newif}%
\@innerdef\@innernewwrite{newwrite}%
\def\@gobble#1{}%
\ifx\inputlineno\@undefined
   \let\@linenumber = \empty % Pre-3.0.
\else
   \def\@linenumber{\the\inputlineno:\space}%
\fi
\def\@futurenonspacelet#1{\def\cs{#1}%
   \afterassignment\@stepone\let\@nexttoken=
}%
\begingroup % The grouping here avoids stepping on an outside use of `\\'.
\def\\{\global\let\@stoken= }%
\\ % now \@stoken is a space token (\\ is a control symbol, so that
\endgroup
\def\@stepone{\expandafter\futurelet\cs\@steptwo}%
\def\@steptwo{\expandafter\ifx\cs\@stoken\let\@@next=\@stepthree
   \else\let\@@next=\@nexttoken\fi \@@next}%
\def\@stepthree{\afterassignment\@stepone\let\@@next= }%
\def\@getoptionalarg#1{%
   \let\@optionaltemp = #1%
   \let\@optionalnext = \relax
   \@futurenonspacelet\@optionalnext\@bracketcheck
}%
\def\@bracketcheck{%
   \ifx [\@optionalnext
      \expandafter\@@getoptionalarg
   \else
      \let\@optionalarg = \empty
      \expandafter\@optionaltemp
   \fi
}%
\def\@@getoptionalarg[#1]{%
   \def\@optionalarg{#1}%
   \@optionaltemp
}%
\def\@nnil{\@nil}%
\def\@fornoop#1\@@#2#3{}%
\def\@for#1:=#2\do#3{%
   \edef\@fortmp{#2}%
   \ifx\@fortmp\empty \else
      \expandafter\@forloop#2,\@nil,\@nil\@@#1{#3}%
   \fi
}%
\def\@forloop#1,#2,#3\@@#4#5{\def#4{#1}\ifx #4\@nnil \else
       #5\def#4{#2}\ifx #4\@nnil \else#5\@iforloop #3\@@#4{#5}\fi\fi
}%
\def\@iforloop#1,#2\@@#3#4{\def#3{#1}\ifx #3\@nnil
       \let\@nextwhile=\@fornoop \else
      #4\relax\let\@nextwhile=\@iforloop\fi\@nextwhile#2\@@#3{#4}%
}%
\@innernewif\if@fileexists
\def\@testfileexistence{\@getoptionalarg\@finishtestfileexistence}%
\def\@finishtestfileexistence#1{%
   \begingroup
      \def\extension{#1}%
      \immediate\openin0 =
         \ifx\@optionalarg\empty\jobname\else\@optionalarg\fi
         \ifx\extension\empty \else .#1\fi
         \space
      \ifeof 0
         \global\@fileexistsfalse
      \else
         \global\@fileexiststrue
      \fi
      \immediate\closein0
   \endgroup
}%
\toks0 = {%
\def\bibliographystyle#1{%
   \@readauxfile
   \@writeaux{\string\bibstyle{#1}}%
}%
\let\bibstyle = \@gobble
\let\bblfilebasename = \jobname
\def\bibliography#1{%
   \@readauxfile
   \@writeaux{\string\bibdata{#1}}%
   \@testfileexistence[\bblfilebasename]{bbl}%
   \if@fileexists
      \nobreak
      \@readbblfile
   \fi
}%
\let\bibdata = \@gobble
\def\nocite#1{%
   \@readauxfile
   \@writeaux{\string\citation{#1}}%
}%
\@innernewif\if@notfirstcitation
\def\cite{\@getoptionalarg\@cite}%
\def\@cite#1{%
   \let\@citenotetext = \@optionalarg
   \printcitestart
   \nocite{#1}%
   \@notfirstcitationfalse
   \@for \@citation :=#1\do
   {%
      \expandafter\@onecitation\@citation\@@
   }%
   \ifx\empty\@citenotetext\else
      \printcitenote{\@citenotetext}%
   \fi
   \printcitefinish
}%
\def\@onecitation#1\@@{%
   \if@notfirstcitation
      \printbetweencitations
   \fi
   \expandafter \ifx \csname\@citelabel{#1}\endcsname \relax
      \if@citewarning
         \message{\@linenumber Undefined citation `#1'.}%
      \fi
      \expandafter\gdef\csname\@citelabel{#1}\endcsname{%
         {\tt
            \escapechar = -1
            \nobreak\hskip0pt
            \expandafter\string\csname#1\endcsname
            \nobreak\hskip0pt
         }%
      }%
   \fi
   \csname\@citelabel{#1}\endcsname
   \@notfirstcitationtrue
}%
\def\@citelabel#1{b@#1}%
\def\@citedef#1#2{\expandafter\gdef\csname\@citelabel{#1}\endcsname{#2}}%
\def\@readbblfile{%
   \ifx\@itemnum\@undefined
      \@innernewcount\@itemnum
   \fi
   \begingroup
      \ifx\begin\@undefined
         \def\begin##1##2{%
            \setbox0 = \hbox{\biblabelcontents{##2}}%
            \biblabelwidth = \wd0
         }%
         \let\end = \@gobble % The arg is `thebibliography' again.
      \fi
      \@itemnum = 0
      \def\bibitem{\@getoptionalarg\@bibitem}%
      \def\@bibitem{%
         \ifx\@optionalarg\empty
            \expandafter\@numberedbibitem
         \else
            \expandafter\@alphabibitem
         \fi
      }%
      \def\@alphabibitem##1{%
         \expandafter \xdef\csname\@citelabel{##1}\endcsname {\@optionalarg}%
         \ifx\biblabelprecontents\@undefined
            \let\biblabelprecontents = \relax
         \fi
         \ifx\biblabelpostcontents\@undefined
            \let\biblabelpostcontents = \hss
         \fi
         \@finishbibitem{##1}%
      }%
      \def\@numberedbibitem##1{%
         \advance\@itemnum by 1
         \expandafter \xdef\csname\@citelabel{##1}\endcsname{\number\@itemnum}%
         \ifx\biblabelprecontents\@undefined
            \let\biblabelprecontents = \hss
         \fi
         \ifx\biblabelpostcontents\@undefined
            \let\biblabelpostcontents = \relax
         \fi
         \@finishbibitem{##1}%
      }%
      \def\@finishbibitem##1{%
         \biblabelprint{\csname\@citelabel{##1}\endcsname}%
         \@writeaux{\string\@citedef{##1}{\csname\@citelabel{##1}\endcsname}}%
         \ignorespaces
      }%
      \let\em = \bblem
      \let\newblock = \bblnewblock
      \let\sc = \bblsc
      \frenchspacing
      \clubpenalty = 4000 \widowpenalty = 4000
      \tolerance = 10000 \hfuzz = .5pt
      \everypar = {\hangindent = \biblabelwidth
                      \advance\hangindent by \biblabelextraspace}%
      \bblrm
      \parskip = 1.5ex plus .5ex minus .5ex
      \biblabelextraspace = .5em
      \bblhook
      \input \bblfilebasename.bbl
   \endgroup
}%
\@innernewdimen\biblabelwidth
\@innernewdimen\biblabelextraspace
\def\biblabelprint#1{%
   \noindent
   \hbox to \biblabelwidth{%
      \biblabelprecontents
      \biblabelcontents{#1}%
      \biblabelpostcontents
   }%
   \kern\biblabelextraspace
}%
\def\biblabelcontents#1{{\bblrm [#1]}}%
\def\bblrm{\rm}%
\def\bblem{\it}%
\def\bblsc{\ifx\@scfont\@undefined
              \font\@scfont = cmcsc10
           \fi
           \@scfont
}%
\def\bblnewblock{\hskip .11em plus .33em minus .07em }%
\let\bblhook = \empty
\def\printcitestart{[}%         left bracket
\def\printcitefinish{]}%        right bracket
\def\printbetweencitations{, }% comma, space
\def\printcitenote#1{, #1}%     comma, space, note (if it exists)
\let\citation = \@gobble
\@innernewcount\@numparams
\def\newcommand#1{%
   \def\@commandname{#1}%
   \@getoptionalarg\@continuenewcommand
}%
\def\@continuenewcommand{%
   \@numparams = \ifx\@optionalarg\empty 0\else\@optionalarg \fi \relax
   \@newcommand
}%
\def\@newcommand#1{%
   \def\@startdef{\expandafter\edef\@commandname}%
   \ifnum\@numparams=0
      \let\@paramdef = \empty
   \else
      \ifnum\@numparams>9
         \errmessage{\the\@numparams\space is too many parameters}%
      \else
         \ifnum\@numparams<0
            \errmessage{\the\@numparams\space is too few parameters}%
         \else
            \edef\@paramdef{%
               \ifcase\@numparams
                  \empty  No arguments.
               \or ####1%
               \or ####1####2%
               \or ####1####2####3%
               \or ####1####2####3####4%
               \or ####1####2####3####4####5%
               \or ####1####2####3####4####5####6%
               \or ####1####2####3####4####5####6####7%
               \or ####1####2####3####4####5####6####7####8%
               \or ####1####2####3####4####5####6####7####8####9%
               \fi
            }%
         \fi
      \fi
   \fi
   \expandafter\@startdef\@paramdef{#1}%
}%
}%
\ifx\nobibtex\@undefined \the\toks0 \fi
\def\@readauxfile{%
   \if@auxfiledone \else % remember: \@auxfiledonetrue if \noauxfile is defined
      \global\@auxfiledonetrue
      \@testfileexistence{aux}%
      \if@fileexists
         \begingroup
            \endlinechar = -1
            \catcode`@ = 11
            \input \jobname.aux
         \endgroup
      \else
         \message{\@undefinedmessage}%
         \global\@citewarningfalse
      \fi
      \immediate\openout\@auxfile = \jobname.aux
   \fi
}%
\newif\if@auxfiledone
\ifx\noauxfile\@undefined \else \@auxfiledonetrue\fi
\@innernewwrite\@auxfile
\def\@writeaux#1{\ifx\noauxfile\@undefined \write\@auxfile{#1}\fi}%
\ifx\@undefinedmessage\@undefined
   \def\@undefinedmessage{No .aux file; I won't give you warnings about
                          undefined citations.}%
\fi
\@innernewif\if@citewarning
\ifx\noauxfile\@undefined \@citewarningtrue\fi
\catcode`@ = \@oldatcatcode
\let\auxfile = \@auxfile
\let\for = \@for
\let\futurenonspacelet = \@futurenonspacelet
\def\iffileexists{\if@fileexists}%
\let\innerdef = \@innerdef
\let\innernewcount = \@innernewcount
\let\innernewdimen = \@innernewdimen
\let\innernewif = \@innernewif
\let\innernewwrite = \@innernewwrite
\let\linenumber = \@linenumber
\let\readauxfile = \@readauxfile
\let\spacesub = \@spacesub
\let\testfileexistence = \@testfileexistence
\let\writeaux = \@writeaux
\def\innerinnerdef#1{\expandafter\innerdef\csname inner#1\endcsname{#1}}%
\innerinnerdef{newbox}%
\innerinnerdef{newfam}%
\innerinnerdef{newhelp}%
\innerinnerdef{newinsert}%
\innerinnerdef{newlanguage}%
\innerinnerdef{newmuskip}%
\innerinnerdef{newread}%
\innerinnerdef{newskip}%
\innerinnerdef{newtoks}%
\def\immediatewriteaux#1{%
  \ifx\noauxfile\@undefined
    \immediate\write\@auxfile{#1}%
  \fi
}%
\begingroup
   \makeactive\^^M \makeactive\ % No spaces or ^^M's from here on.
\gdef\obeywhitespace{%
\makeactive\^^M\def^^M{\par\futurelet\next\@finishobeyedreturn}%
\makeactive\ \let =\ %
\aftergroup\@removebox%
\futurelet\next\@finishobeywhitespace%
}%
\gdef\@finishobeywhitespace{{%
\ifx\next %
\aftergroup\@obeywhitespaceloop%
\else\ifx\next^^M%
\aftergroup\gobble%
\fi\fi}}%
\gdef\@finishobeyedreturn{%
\ifx\next^^M\vskip\blanklineskipamount\fi%
\indent%
}%
\endgroup
\def\@obeywhitespaceloop#1{\futurelet\next\@finishobeywhitespace}%
\def\@removebox{%
  \ifhmode
    \setbox0 = \lastbox
    \ifdim\wd0=\parindent
      \setbox2 = \hbox{\unhbox0}%
      \ifdim\wd2=0pt
        \ignorespaces
      \else
        \box2 % Put it back: it wasn't empty.
      \fi
    \else
       \box0 % Put it back: it wasn't the right width.
    \fi
  \fi
}%
\newskip\blanklineskipamount
\blanklineskipamount = 0pt
\def\frac#1/#2{\leavevmode
   \kern.1em \raise .5ex \hbox{\the\scriptfont0 #1}%
   \kern-.1em $/$%
   \kern-.15em \lower .25ex \hbox{\the\scriptfont0 #2}%
}%
\newdimen\hruledefaultheight  \hruledefaultheight = 0.4pt
\newdimen\hruledefaultdepth   \hruledefaultdepth = 0.0pt
\newdimen\vruledefaultwidth   \vruledefaultwidth = 0.4pt
\def\ehrule{\hrule height\hruledefaultheight depth\hruledefaultdepth}%
\def\evrule{\vrule width\vruledefaultwidth}%
\ifx\sc\undefined
    \def\sc{%
      \expandafter\ifx\the\scriptfont\fam\nullfont
        \font\temp = cmr7 \temp
      \else
        \the\scriptfont\fam
      \fi
      \def\uppercasesc{\char\uccode`}%
    }%
\fi
\ifx\uppercasesc\undefined
  \let\uppercasesc = \relax
\fi
\def\TeX{T\kern-.1667em\lower.5ex\hbox{E}\kern-.125emX\spacefactor1000 }%
\ifx\AmS\undefined
    \def\AmS{{\the\textfont2 A}\kern-.1667em\lower.5ex\hbox
        {\the\textfont2 M}\kern-.125em{\the\textfont2 S}}
\fi
\ifx\AMS\undefined \let\AMS=\AmS \fi
\ifx\AmSLaTeX\undefined
    \def\AmSLaTeX{\AmS-\LaTeX}
\fi
\ifx\AMSLaTeX\undefined \let\AMSLaTeX=\AmSLaTeX \fi
\ifx\AmSTeX\undefined
    \def\AmSTeX{$\cal A$\kern-.1667em\lower.5ex\hbox{$\cal M$}%
            \kern-.125em$\cal S$-\TeX}%
\fi
\ifx\AMSTEX\undefined \let\AMSTEX=\AmSTeX \fi
\ifx\AMSTeX\undefined \let\AMSTeX=\AmSTeX \fi
\ifx\BibTeX\undefined
    \def\BibTeX{B{\sc \uppercasesc i\kern-.025em \uppercasesc b}\kern-.08em
                \TeX}%
\fi
\ifx\BIBTeX\undefined \let\BIBTeX=\BibTeX \fi
\ifx\BIBTEX\undefined \let\BIBTEX=\BibTeX \fi
\ifx\LAMSTeX\undefined
    \def\LAMSTeX{L\raise.42ex\hbox{\kern-.3em\the\scriptfont2 A}%
                 \kern-.2em\lower.376ex\hbox{\the\textfont2 M}%
                 \kern-.125em {\the\textfont2 S}-\TeX}%
\fi
\ifx\LamSTeX\undefined \let\LamSTeX=\LAMSTeX \fi
\ifx\LAmSTeX\undefined \let\LAmSTeX=\LAMSTeX \fi
\ifx\LaTeX\undefined
    \def\LaTeX{L\kern-.36em\raise.3ex\hbox{\sc \uppercasesc a}\kern-.15em\TeX}%
\fi
\ifx\LATEX\undefined \let\LATEX=\LaTeX \fi
\ifx\MF\undefined
    \ifx\manfnt\undefined
            \font\manfnt=logo10
    \fi
    \ifx\manfntsl\undefined
            \font\manfntsl=logosl10
    \fi
    \def\MF{{\ifdim\fontdimen1\font>0pt \let\manfnt = \manfntsl \fi
      {\manfnt META}\-{\manfnt FONT}}\spacefactor1000 }%
\fi
\ifx\METAFONT\undefined \let\METAFONT=\MF \fi
\ifx\SLITEX\undefined
    \def\SLITEX{S\kern-.065em L\kern-.18em\raise.32ex\hbox{i}\kern-.03em\TeX}%
\fi
\ifx\SLiTeX\undefined \let\SLiTeX=\SLITEX \fi
\ifx\SliTeX\undefined \let\SliTeX=\SLITEX \fi
\ifx\SLITeX\undefined \let\SLITeX=\SLITEX \fi
\edef\path{\the\catcode`@}%
\catcode`@ = 11
\let\@oldatcatcode = \path
\newcount \c@tcode
\newcount \c@unter
\newif \ifspecialpathdelimiters
\begingroup
\catcode `\ = 10
\gdef \passivesp@ce { }%
\catcode `\ = 13\relax%
\gdef\activesp@ce{ }%
\endgroup
\def \discretionaries % <delim> <chars> <delim>
    {\begingroup
        \c@tcodes = 13
        \discr@tionaries
    }%
\def \discr@tionaries #1% <delim>
    {\def \discr@ti@naries ##1#1% <chars> <delim>
         {\endgroup
          \def \discr@ti@n@ries ####1% <char> or <delim>
              {\if   \noexpand ####1\noexpand #1%
                     \let \n@xt = \relax
               \else
                     \catcode `####1 = 13
                     \def ####1{\discretionary
                                  {\char `####1}{}{\char `####1}}%
                     \let \n@xt = \discr@ti@n@ries
               \fi
               \n@xt
              }%
          \def \discr@ti@n@ri@s {\discr@ti@n@ries ##1#1}%
         }%
     \discr@ti@naries
    }%
\let\pathafterhook = \relax
\def \path
    {\ifspecialpathdelimiters
        \begingroup
        \c@tcodes = 12
        \def \endp@th {\endgroup \endgroup \pathafterhook}%
     \else
        \def \endp@th {\endgroup \pathafterhook}%
     \fi
     \p@th
    }%
\def \p@th #1% <delim>
    {\begingroup
        \tt
        \c@tcode = \catcode `#1
        \discr@ti@n@ri@s
        \catcode `\ = \active
        \expandafter \edef \activesp@ce {\passivesp@ce \hbox {}}%
        \catcode `#1 = \c@tcode
        \def \p@@th ##1#1% <chars> <delim>
            {\leavevmode \hbox {}##1%
             \endp@th
            }%
     \p@@th
    }%
\def \c@tcodes {\afterassignment \c@tc@des \c@tcode}%
\def \c@tc@des
    {\c@unter = 0
     \loop
            \ifnum \catcode \c@unter = \c@tcode
            \else
                \catcode \c@unter = \c@tcode
            \fi
     \ifnum \c@unter < 255
            \advance \c@unter by 1
     \repeat
     \catcode `\ = 10
    }%
\catcode `\@ = \@oldatcatcode
\discretionaries |~!@$%^&*()_+`-=#{}[]:";'<>,.?\/|%
\def\blackbox{\vrule height .8ex width .6ex depth -.2ex \relax}% square bullet
\def\makeblankbox#1#2{%
  \ifvoid0
    \errhelp = \@makeblankboxhelp
    \errmessage{Box 0 is void}%
  \fi
  \hbox{\lower\dp0
    \vbox{\hidehrule{#1}{#2}%
      \kern -#1% overlap rules
      \hbox to \wd0{\hidevrule{#1}{#2}%
        \raise\ht0\vbox to #1{}% vrule height
        \lower\dp0\vtop to #1{}% vrule depth
        \hfil\hidevrule{#2}{#1}%
      }%
      \kern-#1\hidehrule{#2}{#1}%
    }%
  }%
}%
\newhelp\@makeblankboxhelp{Assigning to the dimensions of a void^^J%
  box has no effect.  Do `\string\setbox0=\string\null' before you^^J%
  define its dimensions.}%
\def\hidehrule#1#2{\kern-#1\hrule height#1 depth#2 \kern-#2}%
\def\hidevrule#1#2{%
  \kern-#1%
  \dimen@=#1\advance\dimen@ by #2%
  \vrule width\dimen@
  \kern-#2%
}%
\newdimen\boxitspace \boxitspace = 3pt
\long\def\boxit#1{%
  \vbox{%
    \ehrule
    \hbox{%
      \evrule
      \kern\boxitspace
      \vbox{\kern\boxitspace \parindent = 0pt #1\kern\boxitspace}%
      \kern\boxitspace
      \evrule
    }%
    \ehrule
  }%
}%
\def\numbername#1{\ifcase#1%
   zero%
   \or one%
   \or two%
   \or three%
   \or four%
   \or five%
   \or six%
   \or seven%
   \or eight%
   \or nine%
   \or ten%
   \or #1%
   \fi
}%
\let\@plainnewif = \newif
\let\@plainnewdimen = \newdimen
\let\newif = \innernewif
\let\newdimen = \innernewdimen
\edef\@eplainoldandcode{\the\catcode`& }%
\catcode`& = 11
\toks0 = {%
\edef\thinlines{\the\catcode`@ }%
\catcode`@ = 11
\let\@oldatcatcode = \thinlines
\def\smash@@{\relax % \relax, in case this comes first in \halign
  \ifmmode\def\next{\mathpalette\mathsm@sh}\else\let\next\makesm@sh
  \fi\next}
\def\makesm@sh#1{\setbox\z@\hbox{#1}\finsm@sh}
\def\mathsm@sh#1#2{\setbox\z@\hbox{$\m@th#1{#2}$}\finsm@sh}
\def\finsm@sh{\ht\z@\z@ \dp\z@\z@ \box\z@}
\edef\@oldandcatcode{\the\catcode`& }%
\catcode`& = 11
\def\&whilenoop#1{}%
\def\&whiledim#1\do #2{\ifdim #1\relax#2\&iwhiledim{#1\relax#2}\fi}%
\def\&iwhiledim#1{\ifdim #1\let\&nextwhile=\&iwhiledim
        \else\let\&nextwhile=\&whilenoop\fi\&nextwhile{#1}}%
\newif\if&negarg
\newdimen\&wholewidth
\newdimen\&halfwidth
\font\tenln=line10
\def\thinlines{\let\&linefnt\tenln \let\&circlefnt\tencirc
  \&wholewidth\fontdimen8\tenln \&halfwidth .5\&wholewidth}%
\def\thicklines{\let\&linefnt\tenlnw \let\&circlefnt\tencircw
  \&wholewidth\fontdimen8\tenlnw \&halfwidth .5\&wholewidth}%
\def\drawline(#1,#2)#3{\&xarg #1\relax \&yarg #2\relax \&linelen=#3\relax
  \ifnum\&xarg =0 \&vline \else \ifnum\&yarg =0 \&hline \else \&sline\fi\fi}%
\def\&sline{\leavevmode
  \ifnum\&xarg< 0 \&negargtrue \&xarg -\&xarg \&yyarg -\&yarg
  \else \&negargfalse \&yyarg \&yarg \fi
  \ifnum \&yyarg >0 \&tempcnta\&yyarg \else \&tempcnta -\&yyarg \fi
  \ifnum\&tempcnta>6 \&badlinearg \&yyarg0 \fi
  \ifnum\&xarg>6 \&badlinearg \&xarg1 \fi
  \setbox\&linechar\hbox{\&linefnt\&getlinechar(\&xarg,\&yyarg)}%
  \ifnum \&yyarg >0 \let\&upordown\raise \&clnht\z@
  \else\let\&upordown\lower \&clnht \ht\&linechar\fi
  \&clnwd=\wd\&linechar
  \&whiledim \&clnwd <\&linelen \do {%
    \&upordown\&clnht\copy\&linechar
    \advance\&clnht \ht\&linechar
    \advance\&clnwd \wd\&linechar
  }%
  \advance\&clnht -\ht\&linechar
  \advance\&clnwd -\wd\&linechar
  \&tempdima\&linelen\advance\&tempdima -\&clnwd
  \&tempdimb\&tempdima\advance\&tempdimb -\wd\&linechar
  \hskip\&tempdimb \multiply\&tempdima \@m
  \&tempcnta \&tempdima \&tempdima \wd\&linechar \divide\&tempcnta \&tempdima
  \&tempdima \ht\&linechar \multiply\&tempdima \&tempcnta
  \divide\&tempdima \@m
  \advance\&clnht \&tempdima
  \ifdim \&linelen <\wd\&linechar \hskip \wd\&linechar
  \else\&upordown\&clnht\copy\&linechar\fi}%
\def\&hline{\vrule height \&halfwidth depth \&halfwidth width \&linelen}%
\def\&getlinechar(#1,#2){\&tempcnta#1\relax\multiply\&tempcnta 8
  \advance\&tempcnta -9 \ifnum #2>0 \advance\&tempcnta #2\relax\else
  \advance\&tempcnta -#2\relax\advance\&tempcnta 64 \fi
  \char\&tempcnta}%
\def\drawvector(#1,#2)#3{\&xarg #1\relax \&yarg #2\relax
  \&tempcnta \ifnum\&xarg<0 -\&xarg\else\&xarg\fi
  \ifnum\&tempcnta<5\relax \&linelen=#3\relax
    \ifnum\&xarg =0 \&vvector \else \ifnum\&yarg =0 \&hvector
    \else \&svector\fi\fi\else\&badlinearg\fi}%
\def\&hvector{\ifnum\&xarg<0 \rlap{\&linefnt\&getlarrow(1,0)}\fi \&hline
  \ifnum\&xarg>0 \llap{\&linefnt\&getrarrow(1,0)}\fi}%
\def\&vvector{\ifnum \&yarg <0 \&downvector \else \&upvector \fi}%
\def\&svector{\&sline
  \&tempcnta\&yarg \ifnum\&tempcnta <0 \&tempcnta=-\&tempcnta\fi
  \ifnum\&tempcnta <5
    \if&negarg\ifnum\&yarg>0                   % 3d quadrant; dp > 0
      \llap{\lower\ht\&linechar\hbox to\&linelen{\&linefnt
        \&getlarrow(\&xarg,\&yyarg)\hss}}\else % 4th quadrant; ht > 0
      \llap{\hbox to\&linelen{\&linefnt\&getlarrow(\&xarg,\&yyarg)\hss}}\fi
    \else\ifnum\&yarg>0                        % 1st quadrant; ht > 0
      \&tempdima\&linelen \multiply\&tempdima\&yarg
      \divide\&tempdima\&xarg \advance\&tempdima-\ht\&linechar
      \raise\&tempdima\llap{\&linefnt\&getrarrow(\&xarg,\&yyarg)}\else
      \&tempdima\&linelen \multiply\&tempdima-\&yarg % 2d quadrant; dp > 0
      \divide\&tempdima\&xarg
      \lower\&tempdima\llap{\&linefnt\&getrarrow(\&xarg,\&yyarg)}\fi\fi
  \else\&badlinearg\fi}%
\def\&getlarrow(#1,#2){\ifnum #2 =\z@ \&tempcnta='33\else
\&tempcnta=#1\relax\multiply\&tempcnta \sixt@@n \advance\&tempcnta
-9 \&tempcntb=#2\relax\multiply\&tempcntb \tw@
\ifnum \&tempcntb >0 \advance\&tempcnta \&tempcntb\relax
\else\advance\&tempcnta -\&tempcntb\advance\&tempcnta 64
\fi\fi\char\&tempcnta}%
\def\&getrarrow(#1,#2){\&tempcntb=#2\relax
\ifnum\&tempcntb < 0 \&tempcntb=-\&tempcntb\relax\fi
\ifcase \&tempcntb\relax \&tempcnta='55 \or
\ifnum #1<3 \&tempcnta=#1\relax\multiply\&tempcnta
24 \advance\&tempcnta -6 \else \ifnum #1=3 \&tempcnta=49
\else\&tempcnta=58 \fi\fi\or
\ifnum #1<3 \&tempcnta=#1\relax\multiply\&tempcnta
24 \advance\&tempcnta -3 \else \&tempcnta=51\fi\or
\&tempcnta=#1\relax\multiply\&tempcnta
\sixt@@n \advance\&tempcnta -\tw@ \else
\&tempcnta=#1\relax\multiply\&tempcnta
\sixt@@n \advance\&tempcnta 7 \fi\ifnum #2<0 \advance\&tempcnta 64 \fi
\char\&tempcnta}%
\def\&vline{\ifnum \&yarg <0 \&downline \else \&upline\fi}%
\def\&upline{\hbox to \z@{\hskip -\&halfwidth \vrule width \&wholewidth
   height \&linelen depth \z@\hss}}%
\def\&downline{\hbox to \z@{\hskip -\&halfwidth \vrule width \&wholewidth
   height \z@ depth \&linelen \hss}}%
\def\&upvector{\&upline\setbox\&tempboxa\hbox{\&linefnt\char'66}\raise
     \&linelen \hbox to\z@{\lower \ht\&tempboxa\box\&tempboxa\hss}}%
\def\&downvector{\&downline\lower \&linelen
      \hbox to \z@{\&linefnt\char'77\hss}}%
\def\&badlinearg{\errmessage{Bad \string\arrow\space argument.}}%
\thinlines
\countdef\&xarg     0
\countdef\&yarg     2
\countdef\&yyarg    4
\countdef\&tempcnta 6
\countdef\&tempcntb 8
\dimendef\&linelen  0
\dimendef\&clnwd    2
\dimendef\&clnht    4
\dimendef\&tempdima 6
\dimendef\&tempdimb 8
\chardef\@arrbox    0
\chardef\&linechar  2
\chardef\&tempboxa  2           % \&linechar and \&tempboxa don't interfere.
\let\lft^%
\let\rt_% distinguish between \rt and \lft
\newif\if@pslope % test for positive slope
\def\@findslope(#1,#2){\ifnum#1>0
  \ifnum#2>0 \@pslopetrue \else\@pslopefalse\fi \else
  \ifnum#2>0 \@pslopefalse \else\@pslopetrue\fi\fi}%
\def\generalsmap(#1,#2){\getm@rphposn(#1,#2)\plnmorph\futurelet\next\addm@rph}%
\def\sline(#1,#2){\setbox\@arrbox=\hbox{\drawline(#1,#2){\sarrowlength}}%
  \@findslope(#1,#2)\d@@blearrfalse\generalsmap(#1,#2)}%
\def\arrow(#1,#2){\setbox\@arrbox=\hbox{\drawvector(#1,#2){\sarrowlength}}%
  \@findslope(#1,#2)\d@@blearrfalse\generalsmap(#1,#2)}%
\newif\ifd@@blearr
\def\bisline(#1,#2){\@findslope(#1,#2)%
  \if@pslope \let\@upordown\raise \else \let\@upordown\lower\fi
  \getch@nnel(#1,#2)\setbox\@arrbox=\hbox{\@upordown\@vchannel
    \rlap{\drawline(#1,#2){\sarrowlength}}%
      \hskip\@hchannel\hbox{\drawline(#1,#2){\sarrowlength}}}%
  \d@@blearrtrue\generalsmap(#1,#2)}%
\def\biarrow(#1,#2){\@findslope(#1,#2)%
  \if@pslope \let\@upordown\raise \else \let\@upordown\lower\fi
  \getch@nnel(#1,#2)\setbox\@arrbox=\hbox{\@upordown\@vchannel
    \rlap{\drawvector(#1,#2){\sarrowlength}}%
      \hskip\@hchannel\hbox{\drawvector(#1,#2){\sarrowlength}}}%
  \d@@blearrtrue\generalsmap(#1,#2)}%
\def\adjarrow(#1,#2){\@findslope(#1,#2)%
  \if@pslope \let\@upordown\raise \else \let\@upordown\lower\fi
  \getch@nnel(#1,#2)\setbox\@arrbox=\hbox{\@upordown\@vchannel
    \rlap{\drawvector(#1,#2){\sarrowlength}}%
      \hskip\@hchannel\hbox{\drawvector(-#1,-#2){\sarrowlength}}}%
  \d@@blearrtrue\generalsmap(#1,#2)}%
\newif\ifrtm@rph
\def\@shiftmorph#1{\hbox{\setbox0=\hbox{$\scriptstyle#1$}%
  \setbox1=\hbox{\hskip\@hm@rphshift\raise\@vm@rphshift\copy0}%
  \wd1=\wd0 \ht1=\ht0 \dp1=\dp0 \box1}}%
\def\@hm@rphshift{\ifrtm@rph
  \ifdim\hmorphposnrt=\z@\hmorphposn\else\hmorphposnrt\fi \else
  \ifdim\hmorphposnlft=\z@\hmorphposn\else\hmorphposnlft\fi \fi}%
\def\@vm@rphshift{\ifrtm@rph
  \ifdim\vmorphposnrt=\z@\vmorphposn\else\vmorphposnrt\fi \else
  \ifdim\vmorphposnlft=\z@\vmorphposn\else\vmorphposnlft\fi \fi}%
\def\addm@rph{\ifx\next\lft\let\temp=\lftmorph\else
  \ifx\next\rt\let\temp=\rtmorph\else\let\temp\relax\fi\fi \temp}%
\def\plnmorph{\dimen1\wd\@arrbox \ifdim\dimen1<\z@ \dimen1-\dimen1\fi
  \vcenter{\box\@arrbox}}%
\def\lftmorph\lft#1{\rtm@rphfalse \setbox0=\@shiftmorph{#1}%
  \if@pslope \let\@upordown\raise \else \let\@upordown\lower\fi
  \llap{\@upordown\@vmorphdflt\hbox to\dimen1{\hss % \dimen1=\wd\@arrbox
    \llap{\box0}\hss}\hskip\@hmorphdflt}\futurelet\next\addm@rph}%
\def\rtmorph\rt#1{\rtm@rphtrue \setbox0=\@shiftmorph{#1}%
  \if@pslope \let\@upordown\lower \else \let\@upordown\raise\fi
  \llap{\@upordown\@vmorphdflt\hbox to\dimen1{\hss
    \rlap{\box0}\hss}\hskip-\@hmorphdflt}\futurelet\next\addm@rph}%
\def\getm@rphposn(#1,#2){\ifd@@blearr \dimen@\morphdist \advance\dimen@ by
  .5\channelwidth \@getshift(#1,#2){\@hmorphdflt}{\@vmorphdflt}{\dimen@}\else
  \@getshift(#1,#2){\@hmorphdflt}{\@vmorphdflt}{\morphdist}\fi}%
\def\getch@nnel(#1,#2){\ifdim\hchannel=\z@ \ifdim\vchannel=\z@
    \@getshift(#1,#2){\@hchannel}{\@vchannel}{\channelwidth}%
    \else \@hchannel\hchannel \@vchannel\vchannel \fi
  \else \@hchannel\hchannel \@vchannel\vchannel \fi}%
\def\@getshift(#1,#2)#3#4#5{\dimen@ #5\relax
  \&xarg #1\relax \&yarg #2\relax
  \ifnum\&xarg<0 \&xarg -\&xarg \fi
  \ifnum\&yarg<0 \&yarg -\&yarg \fi
  \ifnum\&xarg<\&yarg \&negargtrue \&yyarg\&xarg \&xarg\&yarg \&yarg\&yyarg\fi
  \ifcase\&xarg \or  % There is no case 0
    \ifcase\&yarg    % case 1
      \dimen@i \z@ \dimen@ii \dimen@ \or % case (1,0)
      \dimen@i .7071\dimen@ \dimen@ii .7071\dimen@ \fi \or
    \ifcase\&yarg    % case 2
      \or % case 0,2 wrong
      \dimen@i .4472\dimen@ \dimen@ii .8944\dimen@ \fi \or
    \ifcase\&yarg    % case 3
      \or % case 0,3 wrong
      \dimen@i .3162\dimen@ \dimen@ii .9486\dimen@ \or
      \dimen@i .5547\dimen@ \dimen@ii .8321\dimen@ \fi \or
    \ifcase\&yarg    % case 4
      \or % case 0,2,4 wrong
      \dimen@i .2425\dimen@ \dimen@ii .9701\dimen@ \or\or
      \dimen@i .6\dimen@ \dimen@ii .8\dimen@ \fi \or
    \ifcase\&yarg    % case 5
      \or % case 0,5 wrong
      \dimen@i .1961\dimen@ \dimen@ii .9801\dimen@ \or
      \dimen@i .3714\dimen@ \dimen@ii .9284\dimen@ \or
      \dimen@i .5144\dimen@ \dimen@ii .8575\dimen@ \or
      \dimen@i .6247\dimen@ \dimen@ii .7801\dimen@ \fi \or
    \ifcase\&yarg    % case 6
      \or % case 0,2,3,4,6 wrong
      \dimen@i .1645\dimen@ \dimen@ii .9864\dimen@ \or\or\or\or
      \dimen@i .6402\dimen@ \dimen@ii .7682\dimen@ \fi \fi
  \if&negarg \&tempdima\dimen@i \dimen@i\dimen@ii \dimen@ii\&tempdima\fi
  #3\dimen@i\relax #4\dimen@ii\relax }%
\catcode`\&=4  % Back to alignment tab
}%
\catcode`& = 4
\toks2 = {%
\catcode`\&=4  % Back to alignment tab
\def\generalhmap{\futurelet\next\@generalhmap}%
\def\@generalhmap{\ifx\next^ \let\temp\generalhm@rph\else
  \ifx\next_ \let\temp\generalhm@rph\else \let\temp\m@kehmap\fi\fi \temp}%
\def\generalhm@rph#1#2{\ifx#1^
    \toks@=\expandafter{\the\toks@#1{\rtm@rphtrue\@shiftmorph{#2}}}\else
    \toks@=\expandafter{\the\toks@#1{\rtm@rphfalse\@shiftmorph{#2}}}\fi
  \generalhmap}%
\def\m@kehmap{\mathrel{\smash@@{\the\toks@}}}%
\def\mapright{\toks@={\mathop{\vcenter{\smash@@{\drawrightarrow}}}\limits}%
  \generalhmap}%
\def\mapleft{\toks@={\mathop{\vcenter{\smash@@{\drawleftarrow}}}\limits}%
  \generalhmap}%
\def\bimapright{\toks@={\mathop{\vcenter{\smash@@{\drawbirightarrow}}}\limits}%
  \generalhmap}%
\def\bimapleft{\toks@={\mathop{\vcenter{\smash@@{\drawbileftarrow}}}\limits}%
  \generalhmap}%
\def\adjmapright{\toks@={\mathop{\vcenter{\smash@@{\drawadjrightarrow}}}\limits}%
  \generalhmap}%
\def\adjmapleft{\toks@={\mathop{\vcenter{\smash@@{\drawadjleftarrow}}}\limits}%
  \generalhmap}%
\def\hline{\toks@={\mathop{\vcenter{\smash@@{\drawhline}}}\limits}%
  \generalhmap}%
\def\bihline{\toks@={\mathop{\vcenter{\smash@@{\drawbihline}}}\limits}%
  \generalhmap}%
\def\drawrightarrow{\hbox{\drawvector(1,0){\harrowlength}}}%
\def\drawleftarrow{\hbox{\drawvector(-1,0){\harrowlength}}}%
\def\drawbirightarrow{\hbox{\raise.5\channelwidth
  \hbox{\drawvector(1,0){\harrowlength}}\lower.5\channelwidth
  \llap{\drawvector(1,0){\harrowlength}}}}%
\def\drawbileftarrow{\hbox{\raise.5\channelwidth
  \hbox{\drawvector(-1,0){\harrowlength}}\lower.5\channelwidth
  \llap{\drawvector(-1,0){\harrowlength}}}}%
\def\drawadjrightarrow{\hbox{\raise.5\channelwidth
  \hbox{\drawvector(-1,0){\harrowlength}}\lower.5\channelwidth
  \llap{\drawvector(1,0){\harrowlength}}}}%
\def\drawadjleftarrow{\hbox{\raise.5\channelwidth
  \hbox{\drawvector(1,0){\harrowlength}}\lower.5\channelwidth
  \llap{\drawvector(-1,0){\harrowlength}}}}%
\def\drawhline{\hbox{\drawline(1,0){\harrowlength}}}%
\def\drawbihline{\hbox{\raise.5\channelwidth
  \hbox{\drawline(1,0){\harrowlength}}\lower.5\channelwidth
  \llap{\drawline(1,0){\harrowlength}}}}%
\def\generalvmap{\futurelet\next\@generalvmap}%
\def\@generalvmap{\ifx\next\lft \let\temp\generalvm@rph\else
  \ifx\next\rt \let\temp\generalvm@rph\else \let\temp\m@kevmap\fi\fi \temp}%
\toksdef\toks@@=1
\def\generalvm@rph#1#2{\ifx#1\rt % append
    \toks@=\expandafter{\the\toks@
      \rlap{$\vcenter{\rtm@rphtrue\@shiftmorph{#2}}$}}\else % prepend
    \toks@@={\llap{$\vcenter{\rtm@rphfalse\@shiftmorph{#2}}$}}%
    \toks@=\expandafter\expandafter\expandafter{\expandafter\the\expandafter
      \toks@@ \the\toks@}\fi \generalvmap}%
\def\m@kevmap{\the\toks@}%
\def\mapdown{\toks@={\vcenter{\drawdownarrow}}\generalvmap}%
\def\mapup{\toks@={\vcenter{\drawuparrow}}\generalvmap}%
\def\bimapdown{\toks@={\vcenter{\drawbidownarrow}}\generalvmap}%
\def\bimapup{\toks@={\vcenter{\drawbiuparrow}}\generalvmap}%
\def\adjmapdown{\toks@={\vcenter{\drawadjdownarrow}}\generalvmap}%
\def\adjmapup{\toks@={\vcenter{\drawadjuparrow}}\generalvmap}%
\def\vline{\toks@={\vcenter{\drawvline}}\generalvmap}%
\def\bivline{\toks@={\vcenter{\drawbivline}}\generalvmap}%
\def\drawdownarrow{\hbox to5pt{\hss\drawvector(0,-1){\varrowlength}\hss}}%
\def\drawuparrow{\hbox to5pt{\hss\drawvector(0,1){\varrowlength}\hss}}%
\def\drawbidownarrow{\hbox to5pt{\hss\hbox{\drawvector(0,-1){\varrowlength}}%
  \hskip\channelwidth\hbox{\drawvector(0,-1){\varrowlength}}\hss}}%
\def\drawbiuparrow{\hbox to5pt{\hss\hbox{\drawvector(0,1){\varrowlength}}%
  \hskip\channelwidth\hbox{\drawvector(0,1){\varrowlength}}\hss}}%
\def\drawadjdownarrow{\hbox to5pt{\hss\hbox{\drawvector(0,-1){\varrowlength}}%
  \hskip\channelwidth\lower\varrowlength
  \hbox{\drawvector(0,1){\varrowlength}}\hss}}%
\def\drawadjuparrow{\hbox to5pt{\hss\hbox{\drawvector(0,1){\varrowlength}}%
  \hskip\channelwidth\raise\varrowlength
  \hbox{\drawvector(0,-1){\varrowlength}}\hss}}%
\def\drawvline{\hbox to5pt{\hss\drawline(0,1){\varrowlength}\hss}}%
\def\drawbivline{\hbox to5pt{\hss\hbox{\drawline(0,1){\varrowlength}}%
  \hskip\channelwidth\hbox{\drawline(0,1){\varrowlength}}\hss}}%
\def\commdiag#1{\null\,
  \vcenter{\commdiagbaselines
  \m@th\ialign{\hfil$##$\hfil&&\hfil$\mkern4mu ##$\hfil\crcr
      \mathstrut\crcr\noalign{\kern-\baselineskip}
      #1\crcr\mathstrut\crcr\noalign{\kern-\baselineskip}}}\,}%
\def\commdiagbaselines{\baselineskip15pt \lineskip3pt \lineskiplimit3pt }%
\def\gridcommdiag#1{\null\,
  \vcenter{\offinterlineskip
  \m@th\ialign{&\vbox to\vgrid{\vss
    \hbox to\hgrid{\hss\smash@@{$##$}\hss}}\crcr
      \mathstrut\crcr\noalign{\kern-\vgrid}
      #1\crcr\mathstrut\crcr\noalign{\kern-.5\vgrid}}}\,}%
\newdimen\harrowlength \harrowlength=60pt
\newdimen\varrowlength \varrowlength=.618\harrowlength
\newdimen\sarrowlength \sarrowlength=\harrowlength
\newdimen\hmorphposn \hmorphposn=\z@
\newdimen\vmorphposn \vmorphposn=\z@
\newdimen\morphdist  \morphdist=4pt
\dimendef\@hmorphdflt 0       % These two dimensions are
\dimendef\@vmorphdflt 2       % defined by \getm@rphposn
\newdimen\hmorphposnrt  \hmorphposnrt=\z@
\newdimen\hmorphposnlft \hmorphposnlft=\z@
\newdimen\vmorphposnrt  \vmorphposnrt=\z@
\newdimen\vmorphposnlft \vmorphposnlft=\z@

\newdimen\hgrid \hgrid=15pt
\newdimen\vgrid \vgrid=15pt
\newdimen\hchannel  \hchannel=0pt
\newdimen\vchannel  \vchannel=0pt
\newdimen\channelwidth \channelwidth=3pt
\dimendef\@hchannel 0         % Defined via the
\dimendef\@vchannel 2         % macro \getch@nnel
\catcode`& = \@oldandcatcode
\catcode`@ = \@oldatcatcode
}%
\let\newif = \@plainnewif
\let\newdimen = \@plainnewdimen
\ifx\noarrow\@undefined \the\toks0 \the\toks2 \fi
\catcode`& = \@eplainoldandcode
\def\environment#1{%
   \ifx\@groupname\@undefined\else
      \errhelp = \@unnamedendgrouphelp
      \errmessage{`\@groupname' was not closed by \string\endenvironment}%
   \fi
   \edef\@groupname{#1}%
   \begingroup
      \let\@groupname = \@undefined
}%
\def\endenvironment#1{%
   \endgroup
   \edef\@thearg{#1}%
   \ifx\@groupname\@thearg
   \else
      \ifx\@groupname\@undefined
         \errhelp = \@isolatedendenvironmenthelp
         \errmessage{Isolated \string\endenvironment\space for `#1'}%
      \else
         \errhelp = \@mismatchedenvironmenthelp
         \errmessage{Environment `#1' ended, but `\@groupname' started}%
         \endgroup % Probably a typo in the names.
      \fi
   \fi
   \let\@groupname = \@undefined
}%
\newhelp\@unnamedendgrouphelp{Most likely, you just forgot an^^J%
   \string\endenvironment.  Maybe you should try inserting another^^J%
   \string\endgroup to recover.}%
\newhelp\@isolatedendenvironmenthelp{You ended an environment X, but^^J%
   no \string\environment{X} to start it is anywhere in sight.^^J%
   You might also be at an \string\endenvironment\space that would match^^J%
   a \string\begingroup, i.e., you forgot an \string\endgroup.}%
\newhelp\@mismatchedenvironmenthelp{You started an environment named X, but^^J%
   you ended one named Y.  Maybe you made a typo in one^^J%
   or the other of the names?}%
\newif\ifenvironment
\def\checkenv{\ifenvironment \errhelp = \@interwovenenvhelp
   \errmessage{Interwoven environments}%
   \egroup \fi
}%
\newhelp\@interwovenenvhelp{Perhaps you forgot to end the previous^^J%
   environment? I'm finishing off the current group,^^J%
   hoping that will fix it.}%
\newtoks\previouseverydisplay
\newdimen\leftdisplayindent
\newif\if@leftdisplays
\def\leftdisplays{%
  \if@leftdisplays\else
    \previouseverydisplay = \everydisplay
    \everydisplay = {\the\previouseverydisplay \leftdisplaysetup}%
    \let\@save@maybedisableeqno = \@maybedisableeqno
    \let\@saveeqno = \eqno
    \let\@saveleqno = \leqno
    \let\@saveeqalignno = \eqalignno
    \let\@saveleqalignno = \leqalignno
    \let\@maybedisableeqno = \relax
    \def\eqno{\hfill\textstyle\enspace}%
    \def\leqno{%
      \hfill
      \hbox to0pt\bgroup
        \kern-\displaywidth
        \kern-\displayindent
        $\aftergroup\@leftleqnoend
    }%
    \@redefinealignmentdisplays
    \@leftdisplaystrue
  \fi
}%
\def\centereddisplays{%
  \if@leftdisplays
    \everydisplay = \previouseverydisplay
    \let\@maybedisableeqno = \@save@maybedisableeqno
    \let\eqno = \@saveeqno
    \let\leqno = \@saveleqno
    \let\eqalignno = \@saveeqalignno
    \let\leqalignno = \@saveleqalignno
    \@leftdisplaysfalse
  \fi
}%
\def\leftdisplaysetup{%
  \hbox to\displaywidth\bgroup
    \strut
    \dimen@ = \parindent
      \advance\dimen@ by \leftdisplayindent
      \advance\dimen@ by \leftskip
    \hskip\dimen@
    \advance\displayindent by \dimen@
    \advance\displaywidth by -\parindent
      \advance\displaywidth by -\leftdisplayindent
      \advance\displaywidth by -\leftskip
    $%
    \advance\hsize by -\displayindent
    \aftergroup\@lefteqend
    \displaystyle
}%
\def\@lefteqend{\hfil\egroup$}% gets inserted between the ending $$
\def\@leftleqnoend{\hss \egroup$}%
\def\@redefinealignmentdisplays{%
  \def\displaylines##1{\displ@y
    \vcenter{%
      \halign{\hbox to\displaywidth{$\@lign\displaystyle####\hfil$\hfil}\crcr
              ##1\crcr}}}%
  \def\eqalignno##1{\displ@y
    \let\noalign = \@lefteqalignonoalign
    \vcenter{%
      \halign to\displaywidth{%
         \hfil $\@lign\displaystyle{####}$\tabskip\z@skip
        &$\@lign\displaystyle{{}####}$\hfil\tabskip\centering
        &\llap{$\@lign####$}\tabskip\z@skip\crcr
        ##1\crcr}}}%
  \def\leqalignno##1{\displ@y
    \let\eqno = \relax
    \vcenter{%
      \halign to\displaywidth{%
         \hfil$\@lign\displaystyle{####}$\tabskip\z@skip
        &$\@lign\displaystyle{{}####}$\hfil\tabskip\centering
        &\kern-\displaywidth
         \rlap{\kern-\displayindent $\@lign####$}%
         \tabskip\displaywidth\crcr
        ##1\crcr}}}%
}%
\let\@primitivenoalign = \noalign
\newtoks\@everynoalign
\def\@lefteqalignonoalign#1{%
  \@primitivenoalign{%
    \advance\leftskip by -\parindent
    \advance\leftskip by -\leftdisplayindent
    \parskip = 0pt
    \parindent = 0pt
    \the\@everynoalign
    #1%
  }%
}%
\def\monthname{%
   \ifcase\month
      \or Jan\or Feb\or Mar\or Apr\or May\or Jun%
      \or Jul\or Aug\or Sep\or Oct\or Nov\or Dec%
   \fi
}%
\def\fullmonthname{%
   \ifcase\month
      \or January\or February\or March\or April\or May\or June%
      \or July\or August\or September\or October\or November\or December%
   \fi
}%
\def\timestring{\begingroup
   \count0 = \time
   \divide\count0 by 60
   \count2 = \count0   % The hour, from zero to 23.
   \count4 = \time
   \multiply\count0 by 60
   \advance\count4 by -\count0   % The minute, from zero to 59.
   \ifnum\count4<10
      \toks1 = {0}%
   \else
      \toks1 = {}%
   \fi
   \ifnum\count2<12
      \toks0 = {a.m.}%
   \else
      \toks0 = {p.m.}%
      \advance\count2 by -12
   \fi
   \ifnum\count2=0
      \count2 = 12
   \fi
   \number\count2:\the\toks1 \number\count4 \thinspace \the\toks0
\endgroup}%
\def\today{\the\day\ \fullmonthname\ \the\year}%
\newskip\abovelistskipamount      \abovelistskipamount = .5\baselineskip
  \newcount\abovelistpenalty      \abovelistpenalty    = 10000
  \def\abovelistskip{\vpenalty\abovelistpenalty \vskip\abovelistskipamount}%
\newskip\interitemskipamount      \interitemskipamount = 0pt
  \newcount\belowlistpenalty      \belowlistpenalty    = -50
\newskip\belowlistskipamount      \belowlistskipamount = .5\baselineskip
  \newcount\interitempenalty      \interitempenalty    = 0
  \def\interitemskip{\vpenalty\interitempenalty \vskip\interitemskipamount}%
\newdimen\listleftindent    \listleftindent = 0pt
\newdimen\listrightindent   \listrightindent = 0pt
\let\listmarkerspace = \enspace
\newtoks\everylist
\newdimen\@listindent
\def\beginlist{%
  \abovelistskip
  \@listindent = \parindent
  \advance\@listindent by \listleftindent
  \advance\leftskip by \@listindent
  \advance\rightskip by \listrightindent
  \itemnumber = 1
  \the\everylist
}%
\def\li{\@getoptionalarg\@finli}%
\def\@finli{%
  \ifx\@optionalarg\empty \else
    \expandafter\writeitemxref\expandafter{\@optionalarg}%
  \fi
  \ifnum\itemnumber=1 \else \interitemskip \fi
  \printitem
  \advance\itemnumber by 1
  \advance\itemletter by 1
  \advance\itemromannumeral by 1
  \ignorespaces
}%
\def\writeitemxref#1{\definexref{#1}\marker{item}}%
\def\printitem{%
  \par
  \vskip-\parskip
  \noindent
  \printmarker\marker
}%
\def\printmarker#1{\llap{\marker \enspace}}%
\newcount\numberedlistdepth
\newcount\itemnumber
\newcount\itemletter
\newcount\itemromannumeral
\def\numberedmarker{%
  \ifcase\numberedlistdepth
      (impossible)%
  \or \printitemnumber
  \or \printitemletter
  \or \printitemromannumeral
  \else *%
  \fi
}%
\def\printitemnumber{\number\itemnumber}%
\def\printitemletter{\char\the\itemletter}%
\def\printitemromannumeral{\romannumeral\itemromannumeral}%
\def\numberedprintmarker#1{\llap{#1) \listmarkerspace}}%
\def\numberedlist{\environment{@numbered-list}%
  \advance\numberedlistdepth by 1
  \itemletter = `a
  \itemromannumeral = 1
  \beginlist
  \let\marker = \numberedmarker
  \let\printmarker = \numberedprintmarker
}%

\newcount\unorderedlistdepth
\def\unorderedmarker{%
  \ifcase\unorderedlistdepth
      (impossible)%
  \or \blackbox
  \or ---%
  \else *%
  \fi
}%
\def\unorderedprintmarker#1{\llap{#1\listmarkerspace}}%
\def\unorderedlist{\environment{@unordered-list}%
  \advance\unorderedlistdepth by 1
  \beginlist
  \let\marker = \unorderedmarker
  \let\printmarker = \unorderedprintmarker
}%
\def\listing#1{%
   \par \begingroup
   \@setuplisting
   \setuplistinghook
   \input #1
   \endgroup
}%
\let\setuplistinghook = \relax
\def\linenumberedlisting{%
  \ifx\lineno\undefined \innernewcount\lineno \fi
  \lineno = 0
  \everypar = {\advance\lineno by 1 \printlistinglineno}%
}%
\def\printlistinglineno{\llap{[\the\lineno]\quad}}%
\def\listingfont{\tt}%
\def\@setuplisting{%
   \uncatcodespecials
   \obeywhitespace
   \makeactive\`
   \makeactive\^^I
   \def^^L{\vfill\eject}%
   \parskip = 0pt
   \listingfont
}%
{%
   \makeactive\`
   \gdef`{\relax\lq}% Defeat ligatures.
}%
{%
   \makeactive\^^I
   \gdef^^I{\hskip8\fontdimen2}%
}%
\def\verbatimescapechar#1{%
  \gdef\@makeverbatimescapechar{%
    \@makeverbatimdoubleescape #1%
    \catcode`#1 = 0
  }%
}%
\def\@makeverbatimdoubleescape#1{%
  \catcode`#1 = \other
  \begingroup
    \lccode`\* = `#1%
    \lowercase{\endgroup \ece\def*{*}}%
}%
\verbatimescapechar\|  % initially escapechar is |
\def\verbatim{\begingroup
  \uncatcodespecials
  \makeactive\` % make space character a single space, not stretchable
  \@makeverbatimescapechar
  \tt\obeywhitespace}

\def\definecontentsfile#1{%
  \ece\innernewwrite{#1file}%
  \ece\innernewif{if@#1fileopened}%
  \ece\let{#1filebasename} = \jobname
  \ece\def{open#1file}{\opencontentsfile{#1}}%
  \ece\def{write#1entry}{\writecontentsentry{#1}}%
  \ece\def{writenumbered#1entry}{\writenumberedcontentsentry{#1}}%
  \ece\innernewif{ifrewrite#1file} \csname rewrite#1filetrue\endcsname
  \ece\def{read#1file}{\readcontentsfile{#1}}%
}%
\definecontentsfile{toc}%
\def\opencontentsfile#1{%
  \csname if@#1fileopened\endcsname \else
     \ece{\immediate\openout}{#1file} = \csname #1filebasename\endcsname.#1
     \ece\global{@#1fileopenedtrue}%
  \fi
}%
\def\writecontentsentry#1#2#3{\writenumberedcontentsentry{#1}{#2}{#3}{}}%
\def\writenumberedcontentsentry#1#2#3#4{%
  \csname ifrewrite#1file\endcsname
    \csname open#1file\endcsname
    \toks0 = {\expandafter\noexpand \csname #1#2entry\endcsname}%
    \def\temp{#3}%
    \toks2 = \expandafter{#4}%
    \edef\cs{\the\toks2}%
    \edef\@wr{%
      \write\csname #1file\endcsname{%
        \the\toks0 % the \toc...entry control sequence
        {\sanitize\temp}% the text
        \ifx\empty\cs\else {\sanitize\cs}\fi % A secondary number, or nothing:
        {\noexpand\folio}% the page number
      }%
    }%
    \@wr
  \fi
  \ignorespaces
}%
\def\readcontentsfile#1{%
   \edef\temp{%
     \noexpand\testfileexistence[\csname #1filebasename\endcsname]{#1}%
   }\temp
   \if@fileexists
      \input \csname #1filebasename\endcsname.#1\relax
      \csname ifrewrite#1file\endcsname \csname open#1file\endcsname \endif
   \fi
}%
\def\endif{\fi}%
\let\ifxrefwarning = \iftrue
\def\xrefwarningtrue{\@citewarningtrue \let\ifxrefwarning = \iftrue}%
\def\xrefwarningfalse{\@citewarningfalse \let\ifxrefwarning = \iffalse}%
\begingroup
  \catcode`\_ = 8
  \gdef\xrlabel#1{#1_x}%
\endgroup
\def\xrdef#1{\definexref{#1}{\noexpand\folio}{page}}%
\def\definexref#1#2#3{%
  \edef\temp{#1}%
  \readauxfile
  \edef\@wr{\noexpand\writeaux{\string\@definelabel{\temp}{#2}{#3}}}%
  \@wr
  \ignorespaces
}%
\def\@definelabel#1#2#3{%
  \expandafter\gdef\csname\xrlabel{#1}\endcsname{#2}%
  \global\setproperty{\xrlabel{#1}}{class}{#3}%
}%
\def\xrefn#1{%
  \readauxfile
  \expandafter \ifx\csname\xrlabel{#1}\endcsname\relax
    \if@citewarning
       \message{\linenumber Undefined label `#1'.}%
    \fi
    \expandafter\def\csname\xrlabel{#1}\endcsname{%
      `{\tt
        \escapechar = -1
        \expandafter\string\csname#1\endcsname
      }'%
    }%
  \fi
  \csname\xrlabel{#1}\endcsname % Always produce something.
}%
\let\refn = \xrefn
\def\@maybewarnref{%
  \ifundefined{amsppt.sty}%
  \else
    \message{Warning: amsppt.sty and Eplain both define \string\ref. See
             the Eplain manual.}%
    \let\amsref = \ref
  \fi
  \let\ref = \eplainref
  \ref
}
\let\ref = \@maybewarnref
\def\eplainref{\@generalref{}}%
\def\refs{\@generalref s}%
\def\@generalref#1#2{%
  \readauxfile
  \edef\temp{\getproperty{\xrlabel{#2}}{class}}%
  \expandafter\ifx\csname \temp word\endcsname\relax \else
    \csname \temp word\endcsname
    #1\penalty\@M \
  \fi
  \refn{#2}%
}%
\newcount\eqnumber
\newcount\subeqnumber
\def\eqdefn{\@getoptionalarg\@fineqdefn}%
\def\@fineqdefn#1{%
  \ifx\@optionalarg\empty
    \global\advance\eqnumber by 1
    \def\temp{\eqconstruct{\number\eqnumber}}%
  \else
    \def\temp{\@optionalarg}%
  \fi
  \global\subeqnumber = 0
  \gdef\@currenteqlabel{#1}%
  \toks0 = \expandafter{\@currenteqlabel}%
  \begingroup
    \def\eqrefn{\noexpand\eqrefn}%
    \edef\temp{\noexpand\@eqdefn{\the\toks0}{\temp}}%
    \temp
  \endgroup
}%
\def\eqsubdefn#1{%
  \global\advance\subeqnumber by 1
  \toks0 = {#1}%
  \toks2 = \expandafter{\@currenteqlabel}%
  \begingroup
    \def\eqrefn{\noexpand\eqrefn}%
    \def\eqsubreftext{\noexpand\eqsubreftext}%
    \edef\temp{%
      \noexpand\@eqdefn
        {\the\toks0}%
        {\eqsubreftext{\eqrefn{\the\toks2}}{\the\subeqnumber}}%
    }%
    \temp
  \endgroup
}%
\def\@eqdefn#1#2{%
  \definexref{#1}{#2}{eq}%
  \@definelabel{#1}{#2}{eq}%
}%
\def\eqdef{\@getoptionalarg\@fineqdef}%
\def\@fineqdef{%
  \toks0 = \expandafter{\@optionalarg}%
  \edef\temp{\noexpand\@eqdef{\noexpand\eqdefn[\the\toks0]}}%
  \temp
}%
\def\eqsubdef{\@eqdef\eqsubdefn}%
\def\@eqdef#1#2{%
  #1{#2}% Define the label.
  \@maybedisableeqno
  \eqno \eqref{#2}% Print the text.
  \@mayberestoreeqno
  \ignorespaces
}%
\let\@mayberestoreeqno = \relax
\def\@maybedisableeqno{%
  \ifinner
    \global\let\eqno = \relax
    \global\let\@mayberestoreeqno = \@restoreeqno
  \fi
}%
\let\@primitiveeqno = \eqno
\def\@restoreeqno{%
  \global\let\eqno = \@primitiveeqno
  \global\let\@mayberestoreeqno = \empty
}%
\let\eqrefn = \xrefn
\def\eqref#1{\eqprint{\eqrefn{#1}}}%
\let\eqconstruct = \identity
\def\eqprint#1{(#1)}%
\def\eqsubreftext#1#2{#1.#2}%
\let\extraidxcmdsuffixes = \empty
\outer\def\defineindex#1{%
  \def\@idxprefix{#1}%
  \for\@idxcmd:=,marked,submarked,name%
                \extraidxcmdsuffixes\do
  {%
    \@defineindexcmd\@idxcmd
  }%
  \ece\innernewwrite{@#1indexfile}%
  \ece\innernewif{if@#1indexfileopened}%
}%
\newif\ifsilentindexentry
\def\@defineindexcmd#1{%
  \@defineoneindexcmd{s}{#1}\silentindexentrytrue
  \@defineoneindexcmd{}{#1}\silentindexentryfalse
}%
\def\@defineoneindexcmd#1#2#3{%
  \toks@ = {#3}%
  \edef\temp{%
    \def
      \expandonce\csname#1\@idxprefix dx#2\endcsname % e.g., \idx or \sidxname.
      {\def\noexpand\@idxprefix{\@idxprefix}% define \@idxprefix
       \expandonce\csname @@#1idx#2\endcsname
      }%
    \def
      \expandonce\csname @@#1idx#2\endcsname{% e.g., \@@idx
        \the\toks@
        \noexpand\@idxgetrange\expandonce\csname @#1idx#2\endcsname
      }%
  }%
  \temp
}%
\let\indexfilebasename = \jobname
\def\@idxwrite#1#2{%
  \csname if@\@idxprefix indexfileopened\endcsname \else
    \expandafter\immediate\openout\csname @\@idxprefix indexfile\endcsname =
      \indexfilebasename.\@idxprefix dx
    \expandafter\global\csname @\@idxprefix indexfileopenedtrue\endcsname
  \fi
  \def\temp{#1}%
  \edef\@wr{%
    \expandafter\write\csname @\@idxprefix indexfile\endcsname{%
      \string\indexentry
      {\sanitize\temp}%
      {\noexpand#2}%
    }%
  }%
  \@wr
  \ifindexproofing \insert\@indexproof{\indexproofterm{#1}}\fi
  \hookrun{afterindexterm}%
  \ifsilentindexentry \expandafter\ignorespaces\fi
}%
\newif\ifindexproofing
\newinsert\@indexproof
\dimen\@indexproof = \maxdimen                  % No limit on number of terms.
\count\@indexproof = 0  \skip\@indexproof = 0pt % They take up no space.
\font\indexprooffont = cmtt8
\def\indexproofterm#1{\hbox{\strut \indexprooffont #1}}%
\let\@plainmakeheadline = \makeheadline
\def\makeheadline{%
  \indexproofunbox
  \@plainmakeheadline
}%
\def\indexsetmargins{%
  \ifx\undefined\outsidemargin
    \dimen@ = 1truein
    \advance\dimen@ by \hoffset
    \edef\outsidemargin{\the\dimen@}%
    \let\insidemargin = \outsidemargin
  \fi
}%
\def\indexproofunbox{%
  \ifvoid\@indexproof\else
    \indexsetmargins
    \rlap{%
      \kern\hsize
      \ifodd\pageno \kern\outsidemargin \else \kern\insidemargin \fi
      \vbox to 0pt{\unvbox\@indexproof\vss}%
    }%
  \fi
}%
\def\idxrangebeginword{begin}%
\def\idxbeginrangemark{(}% the corresponding magic char to go in the idx file
\def\idxrangeendword{end}%
\def\idxendrangemark{)}%
\def\idxseecmdword{see}%
\def\idxseealsocmdword{seealso}%
\newif\if@idxsee
\let\@idxseenterm = \relax
\def\idxpagemarkupcmdword{pagemarkup}%
\let\@idxpagemarkup = \relax
\def\@idxgetrange#1{%
  \let\@idxrangestr = \empty
  \let\@afteridxgetrange = #1%
  \@getoptionalarg\@finidxgetopt
}%
\def\@finidxgetopt{%
  \for\@idxarg:=\@optionalarg\do{%
    \expandafter\@idxcheckpagemarkup\@idxarg=,%
    \ifx\@idxarg\idxrangebeginword
      \def\@idxrangestr{\idxencapoperator\idxbeginrangemark}%
    \else
      \ifx\@idxarg\idxrangeendword
        \def\@idxrangestr{\idxencapoperator\idxendrangemark}%
      \else
        \ifx\@idxarg\idxseecmdword
          \def\@idxpagemarkup{indexsee}%
          \@idxseetrue
        \else
          \ifx\@idxarg\idxseealsocmdword
            \def\@idxpagemarkup{indexseealso}%
            \@idxseetrue
          \else
             \ifx\@idxpagemarkup\relax
               \errmessage{Unrecognized index option `\@idxarg'}%
             \fi
          \fi
        \fi
      \fi
    \fi
  }%
  \@afteridxgetrange
}%
\def\@idxcheckpagemarkup#1=#2,{%
  \def\temp{#1}%
  \ifx\temp\idxpagemarkupcmdword
    \if ,#2, % If #2 is empty, complain.
      \errmessage{Missing markup command to `pagemarkup'}%
    \else
      \def\temp##1={##1}%
      \edef\@idxpagemarkup{\temp\string#2}%
    \fi
  \fi
}%
\def\idxsubentryseparator{!}%
\def\idxencapoperator{|}%
\def\idxmaxpagenum{99999}%
\newtoks\@idxmaintoks
\newtoks\@idxsubtoks
\def\@idxtokscollect{%
  \edef\temp{\the\@idxsubtoks}%
  \edef\@indexentry{%
    \the\@idxmaintoks
    \ifx\temp\empty\else \idxsubentryseparator\the\@idxsubtoks \fi
    \@idxrangestr
  }%
  \if@idxsee
    \@idxseefalse % Reset so the next term won't be a `see'.
    \edef\temp{\noexpand\@finidxtokscollect{\idxmaxpagenum}}%
  \else
    \def\temp{\@finfinidxtokscollect\folio}%
  \fi
  \temp
}%
\def\@finidxtokscollect#1#2{%
  \def\@idxseenterm{#2}%
  \@finfinidxtokscollect{#1}%
}%
\def\@finfinidxtokscollect#1{%
  \ifx\@idxpagemarkup\relax \else
    \toks@ = \expandafter{\@indexentry}%
    \edef\@indexentry{\the\toks@ \idxencapoperator \@idxpagemarkup}%
    \let\@idxpagemarkup = \relax
  \fi
  \ifx\@idxseenterm\relax \else
    \toks@ = \expandafter{\@indexentry}%
    \edef\@indexentry{\the\toks@{\sanitize\@idxseenterm}}%
    \let\@idxseenterm = \relax
  \fi
  \expandafter\@idxwrite\expandafter{\@indexentry}{#1}%
}%
\def\@idxcollect#1#2{%
  \@idxmaintoks = {#1}%
  \@idxsubtoks = {#2}%
  \@idxtokscollect
}%
\def\@idx#1{%
  #1% Produce TERM as output.
  \@idxcollect{#1}{}%
}%
\def\@sidx#1{\@idxmaintoks = {#1}\@getoptionalarg\@finsidx}%
\def\@finsidx{%
  \@idxsubtoks = \expandafter{\@optionalarg}%
  \@idxtokscollect
}%
\def\idxsortkeysep{@}% This `@' is catcode 11, but it doesn't matter.
\def\@idxconstructmarked#1#2#3{%
  \toks@ = {#2}% The control sequence.
  \toks2 = {#3}% The term.
  \edef\temp{\the\toks2 \idxsortkeysep \the\toks@{\the\toks2}}%
  #1 = \expandafter{\temp}%
}%
\def\@idxmarked#1#2{%
  #1{#2}% Produce \CS{TERM} as output.
  \@idxconstructmarked\@idxmaintoks{#1}{#2}%
  \@idxsubtoks = {}%
  \@idxtokscollect
}%
\def\@sidxmarked#1#2{%
  \@idxconstructmarked\toks@{#1}{#2}%
  \edef\temp{{\the\toks@}}%
  \expandafter\@sidx\temp
}%
\def\@idxsubmarked#1#2#3{%
  #1 #2{#3}% produce `TERM \CS{SUBTERM} as output.
  \@sidxsubmarked{#1}{#2}{#3}%
}%
\def\@sidxsubmarked#1#2#3{%
  \@idxmaintoks = {#1}%
  \@idxconstructmarked\@idxsubtoks{#2}{#3}%
  \@idxtokscollect
}%
\def\idxnameseparator{, }% as in `Tachikawa, Elizabeth'
\def\@idxcollectname#1#2{%
  \def\temp{#1}%
  \ifx\temp\empty
    \toks@ = {}%
  \else
    \toks@ = {\idxnameseparator #1}%
  \fi
  \toks2 = {#2}%
  \edef\temp{\the\toks2 \the\toks@}%
}%
\def\@idxname#1#2{%
  #1 #2% Separate the names by a space in the output.
  \@idxcollectname{#1}{#2}%
  \expandafter\@idxcollect\expandafter{\temp}{}%
}%
\def\@sidxname#1#2{%
  \@idxcollectname{#1}{#2}%
  \expandafter\@sidx\expandafter{\temp}%
}%
\let\indexfonts = \relax
\def\readindexfile#1{%
  \edef\@idxprefix{#1}%
  \testfileexistence[\indexfilebasename]{\@idxprefix nd}%
  \iffileexists \begingroup
    \ifx\begin\undefined
      \def\begin##1{\@beginindex}%
      \let\end = \@gobble
    \fi
    \input \indexfilebasename.\@idxprefix nd
    \singlecolumn
  \endgroup
  \else
    \message{No index file \indexfilebasename.\@idxprefix nd.}%
  \fi
}%
\def\@beginindex{%
  \let\item = \@indexitem
  \let\subitem = \@indexsubitem
  \let\subsubitem = \@indexsubsubitem
  \indexfonts
  \doublecolumns
  \parindent = 0pt
  \hookrun{beginindex}%
}%

\newskip\aboveindexitemskipamount  \aboveindexitemskipamount = 0pt plus2pt
\def\aboveindexitemskip{\vskip\aboveindexitemskipamount}%
\def\@indexitem{\begingroup
  \@indexitemsetup
  \leftskip = 0pt
  \aboveindexitemskip
  \penalty-100 % Encourage page breaks before items.
  \def\par{\endgraf\endgroup\nobreak}%
}%
\def\@indexsubitem{%
  \@indexitemsetup
  \leftskip = 1em
}%
\def\@indexsubsubitem{%
  \@indexitemsetup
  \leftskip = 2em
}%
\def\@indexitemsetup{%
  \par
  \hangindent = 1em
  \raggedright
  \hyphenpenalty = 10000
  \hookrun{indexitem}%
}%
\defineindex{i}%
\begingroup
  \catcode `\^^M = \active %
  \gdef\flushleft{%
    \def\@endjustifycmd{\@endflushleft}%
    \def\@eoljustifyaction{\null\hfil\break}%
    \let\@firstlinejustifyaction = \relax
    \@startjustify %
  }%
  \gdef\flushright{%
    \def\@endjustifycmd{\@endflushright}%
    \def\@eoljustifyaction{\break\null\hfil}%
    \def\@firstlinejustifyaction{\hfil\null}%
    \@startjustify %
  }%
  \gdef\center{%
    \def\@endjustifycmd{\@endcenter}%
    \def\@eoljustifyaction{\hfil\break\null\hfil}%
    \def\@firstlinejustifyaction{\hfil\null}%
    \@startjustify %
  }%
  \gdef\@startjustify{%
    \parskip = 0pt
    \catcode`\^^M = \active %
    \def^^M{\futurelet\next\@finjustifyreturn}%
    \def\@eateol##1^^M{%
      \def\temp{##1}%
      \@firstlinejustifyaction %
      \ifx\temp\empty\else \temp^^M\fi %
    }%
    \expandafter\aftergroup\@endjustifycmd %
    \checkenv \environmenttrue %
    \par\noindent %
    \@eateol %
  }%
  \gdef\@finjustifyreturn{%
    \@eoljustifyaction %
    \ifx\next^^M%
      \def\par{\endgraf\vskip\blanklineskipamount \global\let\par = \endgraf}%
      \@endjustifycmd %
      \noindent %
      \@firstlinejustifyaction %
    \fi %
  }%
\endgroup
\def\@endflushleft{\unpenalty{\parfillskip = 0pt plus1fil\par}\ignorespaces}%
\def\@endflushright{% Remove the \hfil\null\break we just put on.
   \unskip \setbox0=\lastbox \unpenalty
   {\parfillskip = 0pt \par}\ignorespaces
}%
\def\@endcenter{% Remove the \hfil\null\break we just put on.
   \unskip \setbox0=\lastbox \unpenalty
   {\parfillskip = 0pt plus1fil \par}\ignorespaces
}%
\newcount\abovecolumnspenalty   \abovecolumnspenalty = 10000
\newcount\@linestogo         % Lines remaining to process.
\newcount\@linestogoincolumn % Lines remaining in column.
\newcount\@columndepth       % Number of lines in a column.
\newdimen\@columnwidth       % Width of each column.
\newtoks\crtok  \crtok = {\cr}%
\newcount\currentcolumn
\def\makecolumns#1/#2: {\par \begingroup
   \@columndepth = #1
   \advance\@columndepth by #2
   \advance\@columndepth by -1
   \divide \@columndepth by #2
   \@linestogoincolumn = \@columndepth
   \@linestogo = #1
   \currentcolumn = 1
   \def\@endcolumnactions{%
      \ifnum \@linestogo<2
         \the\crtok \egroup \endgroup \par % End \valign and \makecolumns.
      \else
         \global\advance\@linestogo by -1
         \ifnum\@linestogoincolumn<2
            \global\advance\currentcolumn by 1
            \global\@linestogoincolumn = \@columndepth
            \the\crtok
         \else
            &\global\advance\@linestogoincolumn by -1
         \fi
      \fi
   }%
   \makeactive\^^M
   \letreturn \@endcolumnactions
   \@columnwidth = \hsize
     \advance\@columnwidth by -\parindent
     \divide\@columnwidth by #2
   \penalty\abovecolumnspenalty
   \noindent % It's not a paragraph (usually).
   \valign\bgroup
     &\hbox to \@columnwidth{\strut \hsize = \@columnwidth ##\hfil}\cr
}%
\newcount\footnotenumber
\newdimen\footnotemarkseparation \footnotemarkseparation = .5em
\newskip\interfootnoteskip \interfootnoteskip = 0pt
\newtoks\everyfootnote
\newdimen\footnoterulewidth \footnoterulewidth = 2in
\newdimen\footnoteruleheight \footnoteruleheight = 0.4pt
\newdimen\belowfootnoterulespace \belowfootnoterulespace = 2.6pt
\let\@plainfootnote = \footnote
\let\@plainvfootnote = \vfootnote
\def\vfootnote#1{\insert\footins\bgroup
  \interlinepenalty\interfootnotelinepenalty
  \splittopskip\ht\strutbox % top baseline for broken footnotes
  \advance\splittopskip by \interfootnoteskip
  \splitmaxdepth\dp\strutbox
  \floatingpenalty\@MM
  \leftskip\z@skip \rightskip\z@skip \spaceskip\z@skip \xspaceskip\z@skip
  \everypar = {}%
  \parskip = 0pt % because of the vskip
  \ifnum\@numcolumns > 1 \hsize = \@normalhsize \fi
  \the\everyfootnote
  \vskip\interfootnoteskip
  \indent\llap{#1\kern\footnotemarkseparation}\footstrut\futurelet\next\fo@t
}%
\def\footnoterule{\dimen@ = \footnoteruleheight
  \advance\dimen@ by \belowfootnoterulespace
  \kern-\dimen@
  \hrule width\footnoterulewidth height\footnoteruleheight depth0pt
  \kern\belowfootnoterulespace
  \vskip-\interfootnoteskip
}%
\def\numberedfootnote{%
  \global\advance\footnotenumber by 1
  \@plainfootnote{$^{\number\footnotenumber}$}%
}%
\newdimen\paperheight
\ifnum\mag=1000
  \paperheight = 11in % No magnification (yet).
\else
  \paperheight = 11truein % We already have a magnification.
\fi
\def\topmargin{\afterassignment\@finishtopmargin \dimen@}%
\def\@finishtopmargin{%
  \dimen2 = \voffset            % Remember the old \voffset.
  \voffset = \dimen@ \advance\voffset by -1truein
  \advance\dimen2 by -\voffset  % Compute the change in \voffset.
  \advance\vsize by \dimen2     % Change type area accordingly.
}%
\def\advancetopmargin{%
  \dimen@ = 0pt \afterassignment\@finishadvancetopmargin \advance\dimen@
}%
\def\@finishadvancetopmargin{%
  \advance\voffset by \dimen@
  \advance\vsize by -\dimen@
}%
\def\bottommargin{\afterassignment\@finishbottommargin \dimen@}%
\def\@finishbottommargin{%
  \@computebottommargin         % Result in \dimen2.
  \advance\dimen2 by -\dimen@   % Compute the change in the bottom margin.
  \advance\vsize by \dimen2     % Change the type area.
}%
\def\advancebottommargin{%
  \dimen@ = 0pt \afterassignment\@finishadvancebottommargin \advance\dimen@
}%
\def\@finishadvancebottommargin{%
  \advance\vsize by -\dimen@
}%
\def\@computebottommargin{%
  \dimen2 = \paperheight        % The total paper size.
  \advance\dimen2 by -\vsize    % Less the text size.
  \advance\dimen2 by -\voffset  % Less the offset at the top.
  \advance\dimen2 by -1truein   % Less the default offset.
}%
\newdimen\paperwidth
\ifnum\mag=1000
  \paperwidth = 8.5in % No magnification (yet).
\else
  \paperwidth = 8.5truein % We already have a magnification.
\fi
\def\leftmargin{\afterassignment\@finishleftmargin \dimen@}%
\def\@finishleftmargin{%
  \dimen2 = \hoffset            % Remember the old \hoffset.
  \hoffset = \dimen@ \advance\hoffset by -1truein
  \advance\dimen2 by -\hoffset  % Compute the change in \hoffset.
  \advance\hsize by \dimen2     % Change type area accordingly.
}%
\def\advanceleftmargin{%
  \dimen@ = 0pt \afterassignment\@finishadvanceleftmargin \advance\dimen@
}%
\def\@finishadvanceleftmargin{%
  \advance\hoffset by \dimen@
  \advance\hsize by -\dimen@
}%
\def\rightmargin{\afterassignment\@finishrightmargin \dimen@}%
\def\@finishrightmargin{%
  \@computerightmargin          % Result in \dimen2.
  \advance\dimen2 by -\dimen@   % Compute the change in the right margin.
  \advance\hsize by \dimen2     % Change the type area.
}%
\def\advancerightmargin{%
  \dimen@ = 0pt \afterassignment\@finishadvancerightmargin \advance\dimen@
}%
\def\@finishadvancerightmargin{%
  \advance\hsize by -\dimen@
}%
\def\@computerightmargin{%
  \dimen2 = \paperwidth         % The total paper size.
  \advance\dimen2 by -\hsize    % Less the text size.
  \advance\dimen2 by -\hoffset  % Less the offset at the left.
  \advance\dimen2 by -1truein   % Less the default offset.
}%
\let\@plainm@g = \m@g
\def\m@g{\@plainm@g \paperwidth = 8.5 true in \paperheight = 11 true in}%
\newskip\abovecolumnskip \abovecolumnskip = \bigskipamount
\newskip\belowcolumnskip \belowcolumnskip = \bigskipamount
\newdimen\gutter \gutter = 2pc
\newbox\@partialpage
\newdimen\@columnhsize
\newdimen\@normalhsize
\newdimen\@normalvsize  % The original (before multi-columns) \vsize.
\newtoks\previousoutput
\def\quadcolumns{\@columns4}%
\def\triplecolumns{\@columns3}%
\def\doublecolumns{\@columns2}%
\def\begincolumns#1{\ifcase#1\relax \or \singlecolumn \or \@columns2 \or
                            \@columns3 \or \@columns4 \else \relax \fi}%

\let\@ndcolumns = \relax
\chardef\@numcolumns = 1
\def\@columns#1{%
  \@ndcolumns
  \let\@ndcolumns = \@endcolumns
  \chardef\@numcolumns = #1
  \par                     % Shouldn't start in horizontal mode.
  \previousoutput = \expandafter{\the\output}%
  \@columnhsize = \hsize
  \count@ = \@numcolumns
  \advance\count@ by -1
  \advance\@columnhsize by -\count@\gutter
  \divide\@columnhsize by \@numcolumns
  \output = {\global\setbox\@partialpage =
    \vbox{\unvbox255\vskip\abovecolumnskip}%
  }%
  \pagegoal = \pagetotal
  \eject
  \output = {\@columnoutput}%
  \@normalhsize = \hsize
  \@normalvsize = \vsize
  \hsize = \@columnhsize
  \advance\vsize by -\ht\@partialpage
  \advance\vsize by -\ht\footins
  \ifvoid\footins\else \advance\vsize by -\skip\footins \fi
  \multiply\count\footins by \@numcolumns
  \advance\vsize by -\ht\topins
  \ifvoid\topins\else \advance\vsize by -\skip\topins \fi
  \multiply\count\topins by \@numcolumns
  \global\vsize = \@numcolumns\vsize
}%
\def\gutterbox{\vbox to \dimen0{\vfil\hbox{\hfil}\vfil}}%
\newif\if@forceextraline\@forceextralinefalse
\def\@columnsplit{%
  \splittopskip = \topskip
  \splitmaxdepth = \baselineskip
  \dimen@ = \ht255
    \divide\dimen@ by \@numcolumns
 \if@forceextraline
   \advance\dimen@ by \baselineskip
 \fi
 \begingroup
    \vbadness = 10000
    \global\setbox1 = \vsplit255 to \dimen@  \global\wd1 = \hsize
    \global\setbox3 = \vsplit255 to \dimen@  \global\wd3 = \hsize
    \ifnum\@numcolumns > 2
      \global\setbox5 = \vsplit255 to \dimen@ \global\wd5 = \hsize
    \fi
    \ifnum\@numcolumns > 3
      \global\setbox7 = \vsplit255 to \dimen@ \global\wd7 = \hsize
    \fi
  \endgroup
  \if@forceextraline                         % If this is the first time
  \else                                      % through, save the single
    \setbox\@forcelinebox=\copy\@partialpage % column material.
  \fi
  \setbox0 = \box255
  \global\setbox255 = \vbox{%
    \unvbox\@partialpage
    \ifcase\@numcolumns \relax\or\relax
      \or \hbox to \@normalhsize{\box1\hfil\gutterbox\hfil\box3}%
      \or \hbox to \@normalhsize{\box1\hfil\gutterbox\hfil\box3%
                                      \hfil\gutterbox\hfil\box5}%
      \or \hbox to \@normalhsize{\box1\hfil\gutterbox\hfil\box3%
                                      \hfil\gutterbox\hfil\box5%
                                      \hfil\gutterbox\hfil\box7}%
    \fi
  }%
  \setbox\@partialpage = \box0
}%
\def\@columnoutput{%
  \@columnsplit
  \@recoverclubpenalty
  \hsize = \@normalhsize % Local to \output's group.
  \vsize = \@normalvsize
  \the\previousoutput
  \unvbox\@partialpage
  \penalty\outputpenalty
  \global\vsize = \@numcolumns\@normalvsize
}%
\def\singlecolumn{%
  \@ndcolumns
  \chardef\@numcolumns = 1
  \vskip\belowcolumnskip
  \nointerlineskip
}%
\newbox\@forcelinebox
\def\@endcolumns{%
  \global\let\@ndcolumns = \relax
  \par % Shouldn't start in horizontal mode.
  \global\output = {\global\setbox1 = \box255}%
  \pagegoal = \pagetotal
  \eject                     % Exercise the page builder, i.e., \output.
  \setbox2 = \box1           % Save material in box2 in case of overflow.
  \global\setbox255 = \copy2 % Retrieve what the fake \output set.
  \@columnsplit
  \ifvoid\@partialpage
  \else % There is some left-over.
    \setbox0=\box\@partialpage % Merely to void \@partialpage
    \global\setbox255 = \box2  % Retrieve what the fake \output set.
    \@forceextralinetrue       % Add \forcelinebox to \box255 to save single
    \@columnsplit              % column material.
    \global\setbox255 = \vbox{\box\@forcelinebox\box255}%
  \fi
  \global\vsize = \@normalvsize
  \global\hsize = \@normalhsize
  \global\output = \expandafter{\the\previousoutput}%
  \ifvoid\topins\else\topinsert\unvbox\topins\endinsert\fi
  \unvbox255
}%
\def\@saveclubpenalty{% save the current value of \clubpenalty
  \edef\@recoverclubpenalty{%
     \global\clubpenalty=\the\clubpenalty\relax%
     \global\let\noexpand\@recoverclubpenalty\relax
  }% the \noexpand handles infinite self-reference
}%
\let\@recoverclubpenalty\relax
\newdimen\temp@dimen
\def\columnfill{%
 \par
 \dimen@=\pagetotal   % The height of the text set so far.
 \temp@dimen = \vsize % = \@numcolumns * columnheight
 \divide\temp@dimen by \@numcolumns % find out column height
 \loop
   \ifdim \dimen@ > \temp@dimen
     \advance \dimen@ by -\temp@dimen
     \advance \dimen@ by \topskip % fudge factor
 \repeat
 \advance \temp@dimen by -\dimen@
 \advance \temp@dimen by -\prevdepth
 \@saveclubpenalty
 \clubpenalty=10000\relax
 \hrule height\temp@dimen width0pt depth0pt\relax
 \nointerlineskip
 \par
 \nointerlineskip
 \penalty0\vfil % To allow a page break here.
 \relax
}%
\let\wlog = \@plainwlog
\catcode`@ = \@eplainoldatcode
\def\eplain{t}%
{\edef\plainversion{\fmtversion}%
 \xdef\fmtversion{2.8-beta:  7 May 2000 (and plain \plainversion)}%
}%

\font\bigf=cmb18
\font\eightit=cmti8

\newcount\secnum
\outer\def\section #1\par{
    \global\advance \secnum by 1
    %% \subsecnum=0%
    \proclaimnum=0%
    \vskip0pt plus.3\vsize\penalty-250
    \vskip0pt plus-.3\vsize\bigskip\vskip\parskip
    \message{\the\secnum. #1;}
    \leftline{\bf\the\secnum. #1}
    \nobreak\smallskip\noindent}
\def\slabel#1{\definexref{#1}{\the\secnum}{section}}

\newcount\proclaimnum
\outer\def\proclaim #1. #2\par{%
  \global\advance \proclaimnum by 1%
  \medbreak
  \noindent{\bf #1\enspace\the\secnum.\the\proclaimnum\enspace\sl #2\par}%
  \ifdim\lastskip<\medskipamount
    \removelastskip\penalty55\medskip\fi}
\def\plabel#1{\definexref{#1}{\the\secnum.\the\proclaimnum}{proclaim}}

\let\msbm=\tenmsbm

\def\HH{{\rm{H}}}

\def\OO{{\cal O}}
\def\PP{{\rm I\!P}}
\def\ZZ{\hbox{\msbm Z}}

\def\proof{\par\noindent\begingroup{\it Proof:}\rm\enspace}
\def\endproof{\endgroup\hfill{\rm Q.E.D.}\smallskip}

\topglue 1in
\centerline{\bigf On the quantum cohomology of Fano bundles}
\smallskip
\centerline{\bigf over projective spaces}
\bigskip
\centerline{\bf Vincenzo Ancona and  Marco Maggesi%
  \footnote{$^1$}{\rm
    Both authors were  supported by MURST and GNSAGA-INDAM}}

\vskip 30pt

\rightline{\it Dedicated to G\"unther Trautmann on the occasion of his
60$\,^{\hbox{\eightit th}}\,$birthday.}

\vskip 30pt

\midinsert \narrower
\noindent {\bf Abstract.} In their paper ``{\it Quantum cohomology of
projective bundles over {$\PP^n$}}'' (Trans.\ Am.\ Math.\ Soc.\ (1998)
350:9 3615-3638) Z.~Qin and Y.~Ruan introduce interesting techniques
for the computation of the quantum ring of manifolds which are
projectivized bundles over projective spaces; in particular, in the
case of splitting bundles they prove under some restrictions the
formula of Batyrev about the quantum ring of toric manifolds.  Here we
prove the formula of Batyrev on the quantum Chern-Leray equation for a
large class of splitting bundles.
\endinsert

\bigskip

\section Introduction

Over the last years the study of the quantum cohomology ring of
projective manifolds has become a very useful tool in enumerative
geometry. The theory has by now a satisfactory mathematical foundation
\cite{km,rt}.  In general the actual computation is not easy, and
there is more and more interest in the explicit calculation of
concrete examples.

In particular Qin and Ruan in the paper \cite{qr} develop new
interesting techniques for the computation of the quantum ring of
manifolds which are projective bundles associated to holomorphic
vector bundles on projective spaces; in the case of a splitting vector
bundle the corresponding projective bundles are toric manifolds, so
that the Batyrev formula \cite{b} can be taken into account.

Qin and Ruan prove the Batyrev formula for projective splitting
bundles which are Fano manifolds, under some numerical assumption on
the first Chern class.  Their result has been generalized in
\cite{cmr} to the case of splitting vector bundles over products of
projective spaces.

Let $V = \bigoplus_{i=1}^r \OO(m_i)$ be a vector bundle of rank $r$ on
$\PP^n$ and $X = \PP (V)$ the corresponding projective bundle.  Let
$h$ be the cohomology class of a hyperplane of $\PP^n$ and $\xi$ the
class of the tautological line bundle on $\PP(V)$.  We make no
distinction between $h$ and $\pi^*h$ where $\pi \colon \PP (V) \to
\PP^n$ is the projection.  The cohomology ring $\HH^* (X, \ZZ)$ is
generated by $h$ and $\xi$ with the two relations
\item {(i)} the hyperplane equation
  $$ h^{n+1} = 0; \eqdef{hyperplane}$$
\item {(ii)} the Chern-Leray equation
  $$ \sum_{j=0}^r (-1)^j c_j \xi^{r-j} h^j
       = \prod_{i=1}^r (\xi - m_i h)
       = 0  \eqdef{chern}$$
($c_j=c_j(V)$ for $j=0, \dots, r$ are the Chern classes of $V$).

The formula of Batyrev (as stated by Qin and Ruan) reads as follows:
let $V = \bigoplus_{i=1}^r \OO(m_i)$ where $m_i > 0$ for every $i$.
Then the quantum cohomology ring $Q\HH^*(\PP(V))$ is generated by $h$
and $\xi$ with two relations (to be evaluated in the quantum ring)
$$ \eqalignno{
     h^{n+1} &= \prod_{i=1}^r (\xi - m_i h)^{m_i-1} q_2,
     & \eqdef{bat1}\cr
     \prod_{i=1}^r (\xi-m_i h) &= q_1.
     & \eqdef{bat2}} $$

Qin and Ruan prove the above formulas under the additional assumption
$$ \sum_{i=1}^r m_i <
     \min (2r, (n+1+2r)/2, (2n+2+r)/2).
     \eqdef{addhp} $$
Here we give  a proof of the quantum Chern-Leray equation in the case
of a bundle of the form $V = \OO(1) \oplus \OO(m_2) \oplus \dots
\oplus \OO(m_r)$, where $2 \leq m_2 \dots \leq m_r$.

Though the Batyrev formula has been proved by Givental \cite{g} in a
general context, our direct methods can be of some independent
interest, because they could be applied to non-toric cases.  Moreover
some of our computations work in the non Fano case as well.

The more general problem of computing the quantum ring of any Fano
manifold which is a projectivized bundle over a projective manifold
seems at the moment out of reach; nevertheless the case of rank two
bundles should be affordable, mainly because some of them have been
classified: in \cite{sw} when the base space is a surface, and in
\cite{apw} when the base space is a projective space or a quadric.  In
\cite{am} we are able to perform the computation in the
three-dimensional case.

\section Notations

We consider the projective bundle $\pi \colon \PP(V) \to \PP^n$
associated to a rank-$r$ splitting vector bundle $V = \bigoplus_
{i=1}^r \OO(m_i)$ on $\PP^n$.  We assume that the integers $m_i$
satisfy the inequalities $1 = m_1 \leq m_2 \leq \dots \leq m_r$.  In
$\HH_2 (\PP(V), \ZZ)$ we fix the classes $A_1 = (\xi^{r-2} h^n)_*$ and
$A_2 = (\xi^{r-1} h^{n-1} + (1-c_1) \xi^{r-2} h^n)_*$, where the
symbol $(-)_*$ means the Poincar\'e dual homology class; then $A_1$ is
the class of a line in the fiber of $\pi \colon \PP(V) \to \PP^n$ and
$A_2$ is the class of any rational curve constructed as follows.  Let
$\ell \subset \PP^n$ be a line; a surjective morphism $V|_\ell =
\bigoplus_{i=1}^r \OO_\ell(m_i) \to \OO_\ell(1)$ produces an embedding
$$ \tilde \ell = \PP(\OO_\ell(1)) \hookrightarrow \PP(V|_\ell)
     \subset \PP(V) $$
then $[ \tilde \ell\, ] = A_2 $, indeed $\xi [ \tilde \ell\, ] = 1 =
\xi \cdot A_2$ and $h [ \tilde \ell\, ] = h [\ell\,] = 1 = h [A_2]$.
It is well known that $A_1$ and $A_2$ are the extremal rays of the
Mori cone $\hbox{NE} (\PP(V))$.

The anticanonical divisor of $\PP(V)$ is
$$
-K_{\PP(V)} = r \xi + (n+1 - c_1(V)) h
$$ so that $-K_{\PP(V)} \cdot A_1 = r$ and $-K_{\PP(V)} \cdot A_2 =
n+1+r-c_1$.  The moduli space $\overline M_{0,3}(\PP(V), B)$ of
genus-$0$ stable maps with three marked points of class $B = a A_1 + b
A_2$ has virtual dimension
$$ \eqalign{\hbox{virtdim} \big( \overline M_{0,3}(\PP(V), B) \big)
     &= -K \cdot B + \dim(\PP(V)) \cr
     &= a r + b (n+1+r-c_1)+ n + r - 1. \cr} $$

Let $\{T_i\}_{i=0\dots m}$ be a $\ZZ$-basis of $\HH^* (X, \ZZ)$; in
our setting we will take as $\ZZ$-basis the set of monomials
$\xi^lh^m$ with $0 \leq l \leq r-1$ and $0 \leq m \leq n$.  Let $\hat
T_i$ be the cohomology class dual to $T_i$. The quantum product of two
cohomology classes is given by
$$ \alpha * \beta = \sum_{a,b} \sum_i
       I_{aA_1+bA_2}(\alpha,\beta,T_i) \hat T_i q_1^a q_2^b, $$
where $I_{aA_1+bA_2}(\alpha, \beta, T_i)$ are the genus-$0$
Gromov-Witten invariant and $q_1, q_2$ are formal indeterminates.
The degrees of $q_1$ and $q_2$ are defined as $\deg(q_1) = -K_{\PP(V)}
\cdot A_1 = r$ and $\deg(q_2) = -K_{\PP(V)} \cdot A_2 = n+1+r-c_1$
(which allows the quantum product to be homogeneous).

Some of our computation will take place in the quotient ring $Q\HH^*
(\PP(V)) / (q_1)$, that is, we will consider in the quantum product
only the quantum correction coming from the classes $bA_2$ with $b
\geq 1$ and discard all other contributions.  For $ \alpha, \beta \in
Q\HH^*(\PP(V))$ we will write $ \alpha = \beta \pmod {q_1}$ to mean
that $ \alpha$ and $\beta$ are equal in the quotient $Q\HH^* (\PP(V))
/ (q_1)$.  The same for $\alpha = \beta \pmod{q_2}$.

\section On the quantum Chern-Leray equation

\proclaim Theorem.
\plabel{chern-eq}
Let $V := \bigoplus_{i=1}^r \OO_{\PP^n}(m_i)$.  We suppose that $X =
\PP (V)$ is a Fano manifold (i.e., $c_1 \leq n+r$) and $1= m_1 < m_2
\leq m_3 \dots \leq m_r$.  Then Batyrev formula \eqref{bat2}
$$ \prod_{i=1}^r{\kern-.17em\raise1ex\hbox{*}\kern.17em}
     (\xi - m_i h)= q_1 $$
holds.

\proclaim Lemma.
\plabel{contenuto}
Let $Z$ be the hypersurface of $\PP(V)$
$$ \textstyle
 Z = \PP\big(\bigoplus_{i=2}^r\OO(m_i)\,\big)
 \hookrightarrow \PP(V) $$
corresponding to the canonical projection $V \to \bigoplus_{i=2}^r
\OO(m_i)$.  Then for every curve $T$ of class $bA_2$ with $b\geq 1$ we
have $T \cap Z = \emptyset$.

\proof
By induction on $b$, we can suppose $T$ irreducible.  Let us suppose
by contradiction $T \cap Z \ne \emptyset$.  Since $[Z] = \xi - h$ and
$(\xi-h) \cdot T = 0$ we conclude that $T \subset Z$. Let $\xi_Z =
\OO_{\PP(F(-m_2 + 1))}(1)$ be the tautological line bundle on
$\PP(F(-m_2+1))$ where $F(-m_2+1) = \OO(1) \oplus \OO(1 - m_2 + m_3)
\oplus \dots \oplus \OO(1 - m_2 + m_r)$.  Hence
$$ \xi|_Z = \xi_Z + (m_2-1) h; $$
then $0 = (\xi-h) \cdot T = (\xi_Z + (m_2-2) h) \cdot T$ so that
$\xi_Z \cdot T = b(2-m_2) \leq 0$ which  contradicts the ampleness of
$\xi_Z$.
\endproof

\proclaim Lemma.
\plabel{rep}
For $s\le r-1$ and $l\le n$ the class $(\xi-h) \xi^s h^l$ is
represented by a subvariety of $X$ contained in $Z$.

\proof
Let $p = \pi|_Z \colon Z \to \PP^n$ be the projective subbundle of
$\PP(V)$.  The class $h^l (\xi-h)$ is represented by $p^{-1}
(\PP_0^{n-l})$, where $\PP_0^{n-l}$ is a linear subspace of $\PP^n$,
so that the class $(\xi-h) \xi^s h^l$ coincides with the class $\xi^s
p^{-1}(\PP_0^{n-l}) = \xi^s|_Z \cdot p^{-1} (\PP_0^{n-l})$ which is
represented by an effective subvariety of $Z$, because $\xi|_Z$ is
ample and generated by global sections.
\endproof

\proclaim Corollary.
\plabel{pollicino}
For any integer $b \geq 1$ and any triple $\gamma_1, \gamma_2,
\gamma_3$ of cohomology classes in $X$ the following equality holds:
$$ I_{bA_2}\big((\xi-h)\gamma_1, \gamma_2, \gamma_3\big) = 0. $$

\proof
By the linearity of the Gromov-Witten invariants we can suppose
$\gamma_1 = \xi^s h^l$, with $s\le r-1$ and $l\le n$. By
\refn{contenuto} and \refn{rep} the class $(\xi-h) \gamma_1$ is
represented by a subvariety which intersects no curves of class $b
A_2$.
\endproof

\proclaim Lemma.
\plabel{induz}
For $p+q = m+t$ and $p \leq r$, $m \leq r$, the following equality
holds:
$$ \xi^{*p} * h^{*q} - \xi^{*m} * h^{*t} = \xi^{p} h^{q} - \xi^{m}
   h^{t} \pmod{q_1}. $$

\proof
It is enough to show
$$ \xi^{*(l+1)} * h^{*s} - \xi^{*l} * h^{*{(s+1)}}
        = \xi^{l+1}  h^{s} - \xi^{l} h^{s+1} \pmod{q_1} $$
for any $l,s\geq 0, l\le r-1$.  First we suppose $s>0$; let $\gamma_1
:= \xi^{*l} * h^{*s}$, $\gamma_2 := \xi^{*(l+1)} * h^{*(s-1)}$.  By
induction on $l+s$ we obtain
$$ \gamma_1 - \gamma_2 =   \xi^l  h^s - \xi^{l+1}  h^{s-1} \pmod{q_1}. $$
Then
$$ \eqalign{
     \xi^{*(l+1)} * h^{*s} - \xi^{*l} * h^{*{(s+1)}}
       &= (\gamma_1 - \gamma_2)*h  \pmod{q_1} \cr
       &= \xi^{l+1} h^s - \xi^l h^{s+1} + \sum_{b\ge 1}
         \sum_i I_{bA_2}(\xi^l h^s - \xi^{l+1} h^{s-1}, h, T_i)
           \hat{T_i } q_2^b   \pmod{q_1}. \cr} $$
By corollary \refn{pollicino}, all the $ I_{bA_2}(\xi^l h^s -
\xi^{l+1} h^{s-1}, h, T_i) = I_{bA_2}( - (\xi-h) \xi^l h^{s-1}, h,
T_i) $ vanish.  The same computation, interchanging the role of $\xi$
and $h$ works if $l>0$.  The conclusion follows.
\endproof

\bigskip
\begingroup
\par\noindent
{\sl Proof of theorem \refn{chern-eq}:\/}
Since $\deg(q_1) = r$, if $0 \leq i \leq r$, it follows from lemma
\ref{induz} that
$$ \xi^{*(r-i)} * h^{*i} - h^{*r} = \xi^{r-i} h^i - h^r + t_i q_1 $$
for some integers $t_i$.  By lemma 3.5 and 3.7 of \cite{qr}, $I_{A_1}
(\xi, \xi^{r-1}, \xi^{r-1}h^n) = 1$ and $I_{A_1} (\xi^{l_1} h^{m_1},
\xi^{l_2} h^{m_2}, \alpha) = 0$ for $(l_1+l_2) < r$ and any
$\alpha\in\HH^*(\PP(V))$.  So
$$ \xi^{*(r-i)} * h^{*i} = \xi^{r-i}h^i + \left\{\matrix{
        \hfill q_1 \pmod{q_2} & \hbox{if $i=0$} \hfill\cr
        \hfill 0   \pmod{q_2} & \hbox{if $1 \leq i \leq r$}} \right. $$
that is, $t_0=1$ and $t_i=0$ for $1 \leq i \leq r$.
Let $\bar c_V(t)$ be the polynomial
$$ \bar c_V(t) = \sum_{j=0}^r \bar c_j t^j = \prod_{i=1}^r (1-m_it). $$
Then $\bar c_j$ coincides with the Chern class $c_j$ of $V$ for  $0
\leq j \leq \min(n,r)$.  Since $m_1=1$, $\bar c_V(1) = 0$, that is,
$$ \sum_{i=0}^r (-1)^i \bar c_i = 0. $$
It follows
$$ \eqalign{
     \prod_{i=1}^r{\kern-.16em\raise1ex\hbox{*}\kern.16em} (\xi - m_i h)
        &= \sum_{i=0}^r (-1)^i \bar c_i \xi^{*(r-i)}* h^{*i} \cr
        &= \sum_{i=0}^r \Big[(-1)^i \bar c_i \cdot
             (\xi^{r-i} h^i + h^{*r} - h^r + t_i q_1)\Big] \cr
        &= \Big( \sum_{i=0}^r (-1)^i \bar c_i \Big)
             \cdot (h^{*r} - h^r) + q_1 \cr
        &= q_1.} $$
\endproof

\section References \par
\nocite{*}
\bibliography{QChernLeray.bib}
\bibliographystyle{plain}

\bye